\documentclass[12pt,psamsfonts,leqno,oneside,letterpaper]{amsart}
\usepackage[dvips,text={6.5truein,9truein},left=1truein,top=1truein]{geometry}
\usepackage{amssymb,amsmath,amscd,enumerate}
\usepackage[pdftex]{graphicx}
\usepackage{url}

\usepackage[colorlinks,linkcolor=blue,citecolor=blue,pdfstartview=FitH]{hyperref}
\input xy
\xyoption{all}
\SelectTips{cm}{12}

\usepackage{color}

\theoremstyle{definition}
\newtheorem{para}{}[section]
\newtheorem{subpara}{}[para]
\newtheorem{remark}[para]{Remark}
\newtheorem{remarks}[para]{Remarks}
\newtheorem{notation}[para]{Notation}
\newtheorem{convention}[para]{Convention}
\newtheorem{definition}[para]{Definition}
\newtheorem{definitions}[para]{Definitions}

\newcommand\Alternatives{\begin{enumerate}[(i)]}
\newcommand\EndAlternatives{\end{enumerate}}
\newcommand\Conditions{\begin{enumerate}[(1)]}
\newcommand\EndConditions{\end{enumerate}}

\theoremstyle{plain}
\newtheorem{theorem}[para]{Theorem}
\newtheorem{lemma}[para]{Lemma}
\newtheorem{proposition}[para]{Proposition}
\newtheorem{corollary}[para]{Corollary}
\newtheorem{conjecture}[para]{Conjecture}
\newtheorem{claim}[subpara]{}
\newtheorem*{TheoremA}{Theorem A}
\newtheorem*{TheoremB}{Theorem B}
\newtheorem*{simple corollary}{Corollary}
\newtheorem*{simple theorem}{Theorem}
\numberwithin{equation}{para}
\numberwithin{figure}{section}

\newcommand\Number{\begin{para}}
\newcommand\EndNumber{\end{para}}
\newcommand\Definition{\begin{definition}}
\newcommand\EndDefinition{\end{definition}}
\newcommand\Definitions{\begin{definitions}}
\newcommand\EndDefinitions{\end{definitions}}
\newcommand\Theorem{\begin{theorem}}
\newcommand\EndTheorem{\end{theorem}}
\newcommand\Conjecture{\begin{conjecture}}
\newcommand\EndConjecture{\end{conjecture}}
\newcommand\Remark{\begin{remark}}
\newcommand\EndRemark{\end{remark}}
\newcommand\Remarks{\begin{remarks}}
\newcommand\EndRemarks{\end{remarks}}
\newcommand\Convention{\begin{convention}}
\newcommand\EndConvention{\end{convention}}
\newcommand\Notation{\begin{notation}}
\newcommand\EndNotation{\end{notation}}
\newcommand\Lemma{\begin{lemma}}
\newcommand\EndLemma{\end{lemma}}
\newcommand\Proposition{\begin{proposition}}
\newcommand\EndProposition{\end{proposition}}
\newcommand\Corollary{\begin{corollary}}
\newcommand\EndCorollary{\end{corollary}}
\newcommand\Claim{\begin{claim}}
\newcommand\EndClaim{\end{claim}}
\newcommand\Proof{\begin{proof}}
\newcommand\EndProof{\end{proof}}
\newcommand\Equation{\begin{equation}}
\newcommand\EndEquation{\end{equation}}
\newcommand\NoProof{{\hfill$\square$}}
\newcommand\Bullets{\begin{itemize}}
\newcommand\EndBullets{\end{itemize}}

\newcommand\closure{\overline}
\newcommand\trace{\mathop{\rm trace}}
\newcommand\area{\mathop{\rm area}}

\newcommand\RR{{\bf R}}

\newcommand\inter{\mathop{\rm int}}

\newcommand\tM{\widetilde M}
\newcommand\tB{\widetilde B}

\newcommand\tG{\widetilde G}

\newcommand\tQ{\widetilde Q}
\newcommand\tq{\widetilde q}
\newcommand\tf{\widetilde f}
\newcommand\hC{\widehat C}

\newcommand\tGamma{\widetilde \Gamma}
\newcommand\tS{\widetilde S}
\newcommand\h{\widehat}

\newcommand\calf{{\mathcal F}}

\newcommand\call{{\mathcal L}}

\newcommand\calt{{\mathcal T}}

\newcommand\calg{{\mathcal G}}

\newcommand\calx{{\mathcal X}}

\newcommand\calv{{\mathcal V}}

\newcommand\cn{$\ZZ\times\ZZ$-cusp neighborhood}
\newcommand\ZZ{{\mathbb Z}}

\newcommand\NN{{\mathbb N}}
\newcommand\CC{{\mathbb C}}
\newcommand\QQ{{\mathbb Q}}
\newcommand\HH{{\mathbb H}}

\newcommand\cals{{\mathcal S}}

\newcommand\vol{\mathop{\rm vol}}

\newcommand\length{\mathop{{\rm length}}}

\newcommand\rank{\mathop{{\rm rank}}}

\newcommand\isomplus{\mathop{{\rm Isom}_+}}

\newcommand\zzle{{\rm SL}_2}

\title{Small optimal Margulis numbers force upper volume bounds}
\author{Peter Shalen}
\address{Department of Mathematics, Statistics, and Computer Science (M/C 249)\\  University of Illinois at Chicago\\
  851 S. Morgan St.\\
  Chicago, IL 60607-7045} \email{shalen@math.uic.edu}
\thanks{Partially supported by NSF grant DMS-0906155}
\dedicatory{Dedicated to Jos\'e Montesinos on the occasion of his 65th birthday}

\theoremstyle{plain}

\begin{document}
\maketitle
\begin{abstract}
If $\lambda$ is a positive real number strictly less than $\log3$, there is a positive number $V_\lambda$ such that every orientable hyperbolic $3$-manifold of volume greater than $V_\lambda$ admits $\lambda$ as a Margulis number. If $\lambda<(\log3)/2$, such a $V_\lambda$ can be specified explicitly, and is bounded above by 
$$\lambda\bigg(6+\frac{880}{\log3-2\lambda}\log{1\over\log3-2\lambda}\bigg),$$
where $\log$ denotes the natural logarithm. These results imply that for $\lambda<\log3$, an orientable hyperbolic $3$-manifold that does not have $\lambda$ as a Margulis number has a rank-$2$ subgroup of bounded index in its fundamental group, and in particular has a fundamental group of bounded rank. Again, the bounds in these corollaries can be made explicit if $\lambda<(\log3)/2$.
\end{abstract}

\section{Introduction}

Let $M$ be a (complete) orientable hyperbolic $3$-manifold. Up to isometry, we may identify $M$ with
$\HH^3/\Gamma$, where $\Gamma$ is a discrete,
torsion-free subgroup of $\isomplus(\HH^3)$, uniquely determined
up to conjugacy by the hyperbolic structure of $M$, and isomorphic to $\pi_1(M)$. In this setting we have the following definition:

\begin{definition} A {\it
Margulis number} for $M$ (or for $\Gamma$) is a positive real number $\mu$ such that the following condition holds:
\Claim\label{si si si}
If $P$ is a point of $\HH^3$ and $x,y$ are elements of $\Gamma$ such that $\max(d(P,x\cdot P),d(P,y\cdot P))<\mu$,    then $x$ and $y$ commute.
\EndClaim
\end{definition}
Here, and throughout this paper, $d$ denotes  hyperbolic distance on $\HH^3$.

I refer the reader to the introduction to \cite{finiteness} for general background discussion of Margulis numbers, including external references generalizations to the case of higher-dimensional or non-orientable hyperbolic manifolds. In \cite{finiteness} I pointed out that if $\Gamma\cong\pi_1(M)$ is non-abelian then there is an {\it optimal  Margulis number} for $M$, denoted $\mu(M)$, characterized by the property that a given positive number $\mu$ is a Margulis number for $M$ if and only if $\mu\le\mu(M)$. I also discussed the Margulis Lemma, which implies that there is a constant which is a Margulis number for {\it every} orientable hyperbolic $3$-manifold. In this paper I will denote the largest such constant by $\mu_+(3)$.

In \cite{finiteness} I discussed the problem of finding lower bounds for $\mu(M)$, as $M$ varies over a prescribed class of orientable hyperbolic manifolds, and its significance for the problem of classifying finite-volume orientable hyperbolic $3$-manifolds. Results on this problem include Meyerhoff's lower bound of $0.104\ldots$ for $\mu_+(3)$, given in \cite{meyerhoff} (which should be compared with Marc Culler's upper bound of
$0.616\ldots$); the main result of \cite{hakenmarg}, which implies that $\mu(M)\ge0.286$ for any orientable hyperbolic Haken manifold $M$; the main result of \cite{finiteness}, which asserts that, up to isometry, there are at most finitely many orientable $3$-manifolds with $\mu(M)<0.29$; and Corollaries \ref{born free} and \ref{born to facebook} of this paper, which assert that if every subgroup of rank at most $2$ in $\pi_1(M)$ has infinite index---and in particular if $H_1(M;\QQ)$ has rank at least $3$, or if $M$ is closed and $H_1(M;\ZZ_p)$ has rank at least $4$ for some prime $p$---then $\mu(M)\ge\log3=1.09\ldots$. (These corollaries are well-known consequences of known results, and are included in the present paper for completeness. The essential ingredients are, first, the ``$\log3$ Theorem,'' which is deduced from the results of \cite{paradoxical}, \cite{accs}, \cite{agol}, \cite{cg}, \cite{NS}, and \cite{Oh}, and, second, some topological facts proved in \cite{js} and \cite{sw}.)

It is a standard observation that a lower bound for $\mu(M)$ forces a lower bound for the volume of $M$. Indeed, if $\vol M<\infty$, it is easy to deduce that the $\mu$-thick part of $M$ \cite[Chapter D]{bp} is non-empty, and hence that $M$ contains an isometric copy of a ball $B$ of radius $\mu/2$ in $\HH^3$. In particular, the volume of $B$ is a lower bound for $\vol M$. Various refinements of this estimate are also well known.

The theme of the present paper, exemplified by Theorems A and B below and their corollaries, is that for an orientable hyperbolic $3$-manifold $M$, certain {\it upper} bounds on $\mu(M)$ force {\it upper} bounds on the volume of $M$. From this one can deduce that certin upper bounds on $\mu(M)$ also force upper bounds on certain group-theoretical invariants of $\pi_1(M)$ such as its rank.  

Of course, these results can be reinterpreted as saying that suitable lower bounds on the volume of $M$, or on such group-theoretical invariants as the rank of  $\pi_1(M)$, imply certain lower bounds on $\mu(M)$. From this point of view, the paper can be seen as a contribution to a body of results, discussed above, that give lower bounds for $\mu(M)$ under various restrictions on $M$. 
%\redcomment{Decide how much of this to work in: ``consequences of the results of 
%\cite{paradoxical} and generalized in \cite{accs}; the improvement depends on the %Marden Conjecture, which was proved in \cite{agol} and \cite{cg}, and the free case of %the Bers Density Conjecture, which was proved in \cite{NS} and \cite{Oh}.
%This is deduced via a result from \cite{sw} from Corollary below, which gives the same %estimate for $\mu(M)$ under the assumption \ref{born free} will be easily deduced %from Proposition \ref{life is never simple}, which is a slight refinement of \cite[Theorem VI.4.1]{js}, and the ``$\log3$ Theorem,'' which I will quote in a form proved in''}

The following result, which gives a first illustration of how upper bounds on $\mu(M)$ force upper bounds on the volume of $M$, will be proved in the body of the paper as Theorem \ref{abs tract}.

\begin{TheoremA}\label{abs tract intro}
Let $\lambda$ be a positive real number strictly less than $\log3$. Then there is a constant $V_\lambda$ such that for every orientable hyperbolic $3$-manifold $M$ with $\mu(M)<\lambda$ we have  $\vol M\le V_\lambda$ (and in particular $\vol M<\infty$).
\end{TheoremA} 

The following two corollaries to Theorem A will be proved in the body of the paper as 
Corollaries \ref{synecdoche} and \ref{poor thing}.

\begin{simple corollary}\label{synecdoche intro}
Let $\lambda$ be a positive real number strictly less than $\log3$. Then there is a 
there is a natural number $d_\lambda$ such that for every orientable hyperbolic $3$-manifold $M$ with $\mu(M)<\lambda$, the group $\pi_1(M)$ has a rank-$2$ subgroup of index at most $d_\lambda$.
\end{simple corollary}

\begin{simple corollary}\label{poor thing intro}
Let $\lambda$ be a positive real number strictly less than $\log3$. Then there is a 
natural number $k_\lambda$ such that for every orientable hyperbolic $3$-manifold $M$ with $\mu(M)<\lambda$, the group $\pi_1(M)$ has rank at most $k_\lambda$.
\end{simple corollary}

Note that in Theorem A and the two corollaries stated above, no explicit estimate is given for the constants $V_\lambda$, $d_\lambda$ and $k_\lambda$. As I shall now explain, explicit estimates can be obtained if we replace the assumption $\lambda<\log3$ by the stronger assumption that $\lambda<(\log3)/2$.

The following result will be proved in the body of the paper as Corollary \ref{boingo cuckoo}. (It is a corollary to a more technical result, Theorem \ref{what, no soap?}.)

\begin{TheoremB}
Let $\lambda$ be a positive real number strictly less than $(\log3)/2$. Then for every orientable hyperbolic $3$-manifold $M$ with $\mu(M)<\lambda$ we have
$$\vol M< \lambda\bigg(6+\frac{880}{\log3-2\lambda}\log{1\over\log3-2\lambda}\bigg).$$
\end{TheoremB}

To avoid confusion it may be worth pointing out that the right hand side of the inequality in the conclusion of Theorem B is negative if, say, $\lambda<0.1$. Thus in this case the theorem asserts that $\mu(M)$ cannot be less than $\lambda$. However, this is not new information, as Meyerhoff \cite{meyerhoff} has shown that $\mu(3)>0.1$, and indeed his result is used in the proof of Theorem B.

The following two corollaries to Theorem B will be proved in the body of the paper as 
Corollaries \ref{more cuckoo} and \ref{mostly moxie}.

\begin{simple corollary}\label{more cuckoo intro}
Let $\lambda$ be a positive real number strictly less than $(\log3)/2$. Then 
for every orientable hyperbolic $3$-manifold $M$ with $\mu(M)<\lambda$, the group $\pi_1(M)$ has a rank-$2$ subgroup of index at most 
$$\frac{\lambda}{V_0}\bigg(6+\frac{880}{\log3-2\lambda}\log{1\over\log3-2\lambda}\bigg).$$
\end{simple corollary}

\begin{simple corollary}\label{mostly moxie intro}
Let $\lambda$ be a positive real number strictly less than $(\log3)/2$. Then for every hyperbolic $3$-manifold $M$ with 
$\mu(M)<\lambda$, we have
$$\rank\pi_1(M)\le2+\log_2\bigg(\frac{\lambda}{V_0}\bigg(6+\frac{880}{\log3-2\lambda}\log{1\over\log3-2\lambda}\bigg)\bigg).$$
\end{simple corollary}

In my forthcoming paper \cite{cubic}, Theorem \ref{what, no soap?} will be combined with arguments invoking many other results---the $\log3$ theorem, the algebra of
congruence subgroups, and Beukers and Schlickewei's explicit form of Siegel and Mahler's finiteness theorem for solutions to the unit equation in number fields---to prove the following result:

\begin{simple theorem} Let $K$ be any number field, and let $D$ denote its degree. The number of (isometry classes of) closed, non-arithmetic hyperbolic $3$-manifolds which are $\ZZ_6$-homology $3$-spheres, have trace field $K$, and have optimal Margulis number less than $0.183$ is at most
$141\times 2^{24(D+1)}$.
\end{simple theorem}

(That the number of such isometry classes is finite follows from the main result of \cite{finiteness}. It is the explicit bound which is the content of the theorem above.)

The proof of Theorem A occupies Section \ref{abstract section}. The method of proof is to reduce the result to the $\log3$ Theorem (which I mentioned above) using relatively formal arguments based on standard results about algebraic and geometric convergence. The basic strategy is similar to the one used in \cite{finiteness}. 

The proof of Theorem B is rather easily reduced to the case in which $\pi_1(M)$ is a two-generator group. In this case, the proof involves two steps. The first, which is carried out in Section \ref{short section}, consists of showing that an upper bound $\lambda<(\log3)/2$ on $\mu(M)$ forces an explicit upper bound on the minimal length of a non-trivial relation in the generators of $\pi_1(M)$; this step, which is embodied in Proposition \ref{first one}, is a refinement of the elementary packing arguments that are used, for example, in \cite{sw} to give an elementary proof that $(\log3)/2$ is a Margulis number for any hyperbolic $3$-manifold $M$ such that every two-generator subgroup of $\pi_1(M)$ is free. 

The second step involved in proving Theorem B in the case in which $\pi_1(M)$ is a two-generator group is to show that an upper bound on the minimal length of a non-trivial relation in the generators of $\pi_1(M)$ forces an explicit upper bound on the volume of $M$. This bound is given by Proposition \ref{short to bounded}, the proof of which is the goal of Section \ref{short-bounded section}. The proof of Proposition \ref{short to bounded} is a refinement, in the two-generator case, of the argument used by Cooper in \cite{cooper} to show that the volume of a closed hyperbolic $3$-manifold is bounded by $\pi$ times the length of a presentation of its fundamental group. Adapting Cooper's method to giving a volume bound in terms of the length of a single relation rather than the length of an entire presentation involves
some ingredeints are required that did not appear in \cite{cooper}. 
 One of these is an isoperimetric inequality which was proved by Agol and Liu as Lemma 4.4 of their paper \cite{A-L}, where they applied it to prove a result that is somewhat analogous to, but distinct from,  Proposition \ref{short to bounded}. 

The newest ingredient needed for the proof of Proposition \ref{short to bounded} is a deep and technical topological result, Lemma \ref{oh far out}. One of the results needed to prove Lemma \ref{oh far out}, Proposition \ref{wow what a variety of characters}, is a new application of the characteristic submanifold theory that seems to be of particular independent interest.

In Section \ref{concrete section} I assemble the results of Sections \ref{short section} and \ref{short-bounded section} to prove Theorem \ref{what, no soap?} and its various corollaries (including Corollary \ref{boingo cuckoo} which is Theorem B above).

Sections \ref{gen prelim} and \ref{3 prelim} are preliminary sections, devoting to assembling more or less well known results used in the later sections for which convenient references are not easy to find.

I am grateful to Michael Siler, Dick Canary and Marc Culler for  a series of valuable discussions of the material in this paper. Siler pointed out an error in an earlier version of the paper; he called my attention to Lemma 4.4 of Agol and Liu's paper \cite{A-L}; and most importantly he made a suggestion that led me to the realization that $\pi$ could be replaced by $\min(\pi,\lambda)$ in Proposition \ref{short to bounded}, which gives an important improvement in Theorem \ref{what, no soap?}. Canary has painstakingly continued to educate me about algebraic and geometric convergence; he corrected my na\"\i ve ideas about the algebraic convergence arguments involved in the proof of Theorem \ref{abs tract}, and helped me to the correct argument and suitable references. Culler patiently listened to all the details of the material in the paper before they were written down, and clarified the proof of Proposition \ref{uneeda}.   

\section{General preliminaries}\label{gen prelim}

Throughout this paper, $\log x$ will denote the natural logarithm of a positive number $x$, and $\#(X)$ will denote the cardinality of a finite set $X$.

In statements and arguments involving fundamental groups, I will suppress base points
whenever it is possible to do so without ambiguity. If $X$ is a path-connected space, I will often implicitly assume that $X$ is equipped with an unnamed (and arbitrary) base point $\star_X$, and write $\pi_1(X)$ for $\pi_1(X,\star_X)$. If $f:X\to Y$ is a map between path-connected spaces then  $f_\sharp:\pi_1(X)\to\pi_1(Y)$ will mean the homomorphism  from $\pi_1(X,\star_X)$ to $\pi_1(Y)=\pi_1(Y,\star_Y)$ which is obtained by composing the standard induced homomorphism $f_\sharp:\pi_1(X,\star_X)\to\pi_1(Y,f(\star_X))$ with the isomorphism from $\pi_1(Y,f(\star_X))$ to $\pi_1(Y,\star_Y)$ defined by an unspecified path from $f(\star_X)$ to $\star_Y$. Thus $f_\sharp$ is well-defined up to post-composition with inner isomorphisms of $\pi_1(Y)$. Many assertions about $f_\sharp$, such as the assertion that is injective or surjective, are invariant
 under post-composition with inner isomorphisms, and will be made without reference to a connecting path. Likewise, the image of an element or subgroup of $\pi_1(X)$ under $f_\sharp$ is well-defined up to conjugacy.

A path connected subset $A$ of a path connected space $X$ will be termed {\it $\pi_1$-injective} if the inclusion homomorphism $\pi_1(A)\to\pi_1(X)$ is injective.

The following two easy results from group theory will be needed later in the paper.

\Lemma\label{o hula who}
Let $F$ be a free group on a generating set $S$, let $Z$ be a cyclic subgroup of $F$, and let $k$ be a positive integer. Then there are at most $2k+1$ elements of $Z$ that can be expressed as words of length at most $k$ in the generating set $S$.
\EndLemma

\Proof
Let $t$ be a generator of $Z$. If $t=1$ the assertion is trivial. If $t\ne1$, there exist a reduced word $V$ in the generating set $S$ and a cyclically reduced word $W$ in $S$ such that the word $V\ast W\ast\overline V$ is reduced and represents $t$; here $\ast$ denotes concatenation of words and $\overline V$ denotes the inverse of the word $V$. For any non-zero integer $n$ the word $V\ast(\star^n W)\ast\overline V$ is reduced and represents $t^n$; here $\star^nW$ denotes the $n$-fold concatenation of $W$. In particular, the unique reduced word representing $t^n$ has length at least $n$. If $t^n$ can be expressed as a word of length at most $k$ then the reduced word represented $t^n$ has length at most $k$, and hence $n\le k$. The conclusion follows.
\EndProof

\Proposition\label{you peeked}
Suppose that $\tGamma$ is a finite-index subgroup of a finitely generated group $\Gamma$. Then
$$\rank\Gamma\le \rank\tGamma+\log_2[\Gamma:\tGamma].$$
\EndProposition

\Proof
Set $r=\rank\tGamma$, and fix a generating set $\tS$ for $\tGamma$ with $|\tS|=r$. Let $S\supset\tS$ be a finite generating set for $\Gamma$, and let $S$ be chosen to have minimal cardinality among all generating sets for $\Gamma$ that contain $S$. Let us denote the distinct elements of $S-\tS$ by $x_1,\ldots,x_k$. For $0\le j\le k$, let $\Gamma_j$ denote the subgroup of $\Gamma$ generated by $\tS\cup\{x_1,\ldots,x_j\}$ (so that in particular $\Gamma_0=\tGamma$ and $\Gamma_k=\Gamma$). It follows from the minimality of $S$ that $\Gamma_{j-1}$ is a proper subgroup of $\Gamma_j$ for $j=1,\ldots,k$, and therefore $[\Gamma_j:\Gamma_{j-1}]\ge2$. Hence
$$[\Gamma:\tGamma]=\prod_{j=1}^k
[\Gamma_j:\Gamma_{j-1}]\ge2^k.$$
Using this, we find
$$\rank\Gamma\le|S|=r+k\le r+ \log_2[\Gamma:\tGamma].$$
\EndProof

The following elementary fact from hyperbolic geometry will also be needed.

\Proposition\label{uneeda}
Let $T$ be a triangle in a hyperbolic space, and let $L$ denote the length of the shortest side of $T$. Then
$$ \area T<\min(\pi,L).$$
\EndProposition

\Proof
We may assume that $T\subset\HH^2$. let $l$ denote a side of $T$ having length $L$, and let $\overline\HH^2$ denote the union of $\HH^2$ with the circle at infinity. There is a triangle $T'$ in $\overline\HH^2$ which contains $T$, has $l$ as a side, and has an ideal vertex opposite $l$. It is enough to prove that $\area T'<\min(\pi,L)$. Let us identify $\HH^2$ with the upper half-plane model in such a way that $l$ is an arc in the upper unit semicircle and the other two sides of $T'$ are vertical rays. Let $(x_1,y_1)$ and $(x_2,y_2)$ be the endpoints of $l$ in Cartesian coordinates, where $-1<x_1<x_2<1$. We have
$$\area T'=\int_{x_1}^{x_2}\int_{\sqrt{1-x^2}}^\infty\frac1{y^2}\,dy\,dx=\arcsin x_2-\arcsin x_1,$$
which is the Euclidean length of the arc $l$. This shows that $\area T'<\pi$. On the other hand, the hyperbolic length $L$ of $l$ is the integral over the arc $l$ of the hyperbolic length element, which is given by $ds/y$ where $ds$ is the Euclidean length element; since $y<1$ at all but at most one point of the arc $l$, the hyperbolic length of $l$ strictly exceeds its Euclidean length, so that $\area T'<L$.
\EndProof

%\Definitions\label{colonoscopy}
%By a {\it hyperbolic polytope} in $\HH^n$ I will mean a non-empty compact subset of %$\HH^n$ which is a finite intersection of hyperbolic hyperplanes and hyperbolic half-spaces. If $M$ is a hyperbolic $n$-manifold, a {\it hyperbolic polytope} in $M$ will mean a subset $P$ of $M$ such that a local isometry $\HH^n\to M$ maps some hyperbolic polytope in $\HH^n$ homeomorphically onto $P$.

%Let $M$ be a hyperbolic manifold. A {\it hyperbolic PL structure} on $M$ is a PL structure with respect to which every hyperbolic polytope in $M$ is a PL set.

%\EndDefinitions

%\Proposition\label{a polyp case} Every finite-volume hyperbolic manifold admits a hyperbolic PL structure. 
%\EndProposition

%\Proof
%\redcomment{Give it.}
%\EndProof

\section{Three-manifold preliminaries}\label{3 prelim}

\Number\label{cat-a-gory}

When no category is specified, it will be understood that ``manifolds'' and ``submanifolds'' are smooth. At various points in the paper I will need to mention PL or real-analtyic manifolds, and in these cases I will be explicit about the category in question. 

As is customary in $3$-manifold topology, I will frequently quote results about $3$-manifolds that are proved in the PL category, and apply them in the smooth category, most often without explicitly mentioning the transition. In most cases this will be justified by the following facts:
\begin{enumerate}
\item \label{sunny side up}If $M$ is a smooth $n$-manifold and $V$ is a (possibly empty) smooth, properly embedded submanifold, then $M$ has a smooth triangulation with respect to which $V$ is a polyhedral subset.
\item\label{over easy} Any two smooth triangulations of a given smooth manifold determine the same PL structure up to PL homeomorphism.
\item\label{coddled} If $n\le3$, and if $M$ and $M'$ are smooth $n$-manifolds which are PL homeomorphic with respect to the PL structures defined by smooth triangulations, then $M$ and $M'$ are diffeomorphic.
\end{enumerate}
Of these facts, (\ref{over easy}), and the case of (\ref{sunny side up}) where $V=\emptyset$, are proved in \cite{whitehead}, and (\ref{coddled}) is included in \cite[Theorem 6.3]{munkres}. I have not located a direct reference for the case of
(\ref{sunny side up}) in which $V\not=\emptyset$, but it is a very special case of the main result of \cite{goresky}.
\EndNumber

In Section \ref{short-bounded section} I will use a somewhat different result about the interaction between the PL and smooth (or rather real-analytic) categories:

\Proposition \label{real analytic}
Suppose that $M$ is a real-analytic manifold and that $X_1,\ldots,X_n$ are semi-analytic sets. Then $M$ admits a smooth triangulation with respect to which $X_1,\ldots,X_n$ are polyhedral subsets.
\EndProposition

\Proof
It is shown in \cite{grauert} that $M$ is real-analytically isomorphic to a real-analytic submanifold $M'$ of $\RR^N$ for some $N$. If we identify $M$ with $M'$, the sets 
$M,X_1,\ldots,X_n$ become semi-analtyic subsets of $\RR^N$. The main theorem of
of \cite{graver} then asserts that $\RR^N$ admits a smooth triangulation with respect to which $M,X_1,\ldots,X_n$ are polyhedral subsets. The conclusion follows.
\EndProof

\Corollary\label{gruosso}
Let $M$ be a hyperbolic $n$-manifold, let $p:\HH^n\to M$ be a locally isometric covering map, and let $\delta_1,\ldots,\delta_k$ be closed hyperbolic simplices in $\HH^n$. Then
$M$ admits a smooth triangulation with respect to which $X_1,\ldots,X_n$ are polyhedral subsets. 
\EndCorollary

\Proof
Each $\delta_i$ may be subdivided into finitely many closed hyperbolic simplices which are embedded in $M$ under $p$. Hence we may assume without loss of generality that the $\delta_i$ are themselves embedded in $M$ under $p$. In this case the sets $p(\delta_1),\ldots,p(\delta_k)$ are clearly semi-analytic in the real analytic structure defined by the hyperbolic structure of $M$, and so the result follows from Proposition \ref{real analytic}.
\EndProof

I will say that an orientable $3$-manifold $M$ is {\it irreducible} if $M$ is connected and every (smooth) $2$-sphere in $M$ bounds a (smooth) ball in $M$.

\Lemma\label{plug}Let $X_0$ be a compact, connected, $3$-dimensional submanifold of an irreducible, orientable $3$-manifold $M$. Then there is a compact, irreducible, $3$-dimensional submanifold $X_1$ of $M$ such that $X_1\supset X_0$, and such that the inclusion homomorphism $\pi_1(X_0)\to\pi_1(X_1)$ is surjective.  \EndLemma

\Proof Let $\calx$ denote the set of all compact $3$-dimensional submanifolds $X$ of $M$ such that $X\supset X_0$ and the inclusion homomorphism $\pi_1(X_0)\to\pi_1(X)$ is surjective. We have $X_0\in\calx$. Let $X_1\in\calx$ be chosen to have the smallest number of boundary components among all submanifolds in $\calx$. It suffices to show that $X_1$ is irreducible. If $S\subset\inter X_1$ is a $2$-sphere then $S$ bounds a ball $B\subset Q$. If $B\not\subset X_1$, then $X_1\cup B$ is a compact submanifold of $Q$ containing $X_1$ and having fewer boundary components than $X_1$. It is clear that the inclusion homomorphism $\pi_1(X_1)\to\pi_1(X_1\cup B)$ is surjective, and hence that $X_1\cup B\in\calx$, a contradiction to minimality.  \EndProof

\Definitions\label{simple def}
Let $M$ be a, irreducible, orientable $3$-manifold. 
Following the convention of \cite{hempel}, I will define an {\it incompressible surface} in $M$ to be a compact, properly embedded $2$-manifold in $M$ which if it is $\pi_1$-injective and is not a sphere or a boundary-parallel disk.

A closed (smooth) $2$-manifold $S$ in $\int M$ is said to be {\it boundary-parallel} if $S$ is the frontier of a (smooth) submanifold $H$ of $M$ such that the pair $(H,S)$ is diffeomorphic to $(S\times[0,1],S\times\{1\})$. (The definition in the case of a properly embedded surface with non-empty boundary would be slightly trickier in the smooth category, but will not be needed in this paper.)

I will define a {\it Haken manifold} to be a compact, irreducible, orientable $3$-manifold which contains an incompressible surface. Note that according to this definition a $3$-ball is not a Haken manifold. 

By an {\it essential disk} in the irreducible, orientable $3$-manifold $M$ I will mean a properly embedded disk whose boundary does not bound a disk in $\partial M$.  I will say that $M$ is {\it boundary-irreducible} if $M$ contains no essential disk. It follows from the Loop Theorem \cite[p. 39, 4.2]{hempel} that $M$ is boundary-irreducible if and only if every component $T$ of $\partial M$ is $\pi_1$-injective in $M$.  \EndDefinitions

\Lemma\label{it's parallel}
Let $N$ be an irreducible, orientable $3$-manifold, let $T\subset\inter M$ be a closed incompressible surface, and suppose that the inclusion map $T\to N$ is homotopic to a map of $T$ into $\partial N$. Then $T$ is boundary-parallel in $N$.
\EndLemma

\Proof
In the case where $N$ is compact this is included in \cite[Lemma 5.3]{waldhausen}. To prove it in the general case, note that since  the inclusion map $i:T\to N$ is homotopic to some map $f$ of $T$ into $\partial N$, there is a compact subset $X_0$ of $N$, containing $T$ and $f(T)$, such that $i$ is homotopic to $f$ in $X_0$. After passing to a regular neighborhood we may assume that $X_0$ is a submanifold of $N$. Accoding to Lemma \ref{plug}, there is a compact
irreducible submanifold $X_1$ of $M$ such that $X_1\supset X_0$. Applying the compact case of the lemma, with $X_1$ in place of $N$, we deduce that $T$ is boundary parallel in $X_1$, and hence in $N$.
\EndProof

\Definition Let $M$ be an orientable hyperbolic $3$-manifold, and let $p:\HH^3\to M$ be a locally isometric covering map. 
A {\it \cn} in $M$ is a subset $C$ of $M$ such that $p^{-1}(C)\subset\HH^3$ is a horoball, and the image of the inclusion homomorphism $\pi_1(C)\to\pi_1(M)$ is a free abelian group of rank $2$. 
\EndDefinition

\Proposition\label{hypersimple} Let $M$ be a hyperbolic $3$-manifold. Then 
\begin{enumerate}
\item\label{it's irreducible}
$M$ is aspherical and irreducible. 
%\item\ref{cyclic knot}
%For every non-cyclic abelian subgroup $X$ of $\pi_1(M)$, there is a  cusp %neighborhood $C\subset M$ such that $X$ is conjugate to a subgroup of the image of %the inclusion homomorphism $\pi_1(C)\to\pi_1(M)$. 
\item\label{yes injective}
Every incompressible torus $T\subset M$ is the boundary of a
submanifold $C$ of $M$ which is closed as a subset of $M$, is diffeomorphic to $T^2\times[0,\infty)$, and has finite volume. 
\item\label{not injective}
For every torus $T\subset M$ which is {\it not} incompressible, either
(a) $T$ is contained in a  $3$-ball in $M$, or (b) $T$ is the boundary of a solid torus in $M$.
\end{enumerate}
\EndProposition

\Proof
Write $M=\HH^3/\Gamma$ where $\Gamma\subset\isomplus( \HH^3)$ is discrete and torsion-free. Let $p:\HH^3\to M$ denote the quotient map.
Since the universal covering $\HH^3$ of $M$ is contractible, $M$ is aspherical. 
To prove that $M$ is irreducible, suppose that $S\subset M$ is a $2$-sphere. Then $S$ lifts to a $2$-sphere $\tS\subset\HH^3$. Since $\HH^3$ is diffeomorphic to $\RR^3$, it is irreducible by \cite[Theorem 1]{moise}, and so $\tS$ bounds a ball $\tB\subset\HH^3$. For any $\gamma\in\Gamma-\{1\}$, we have $\tS\cap\gamma\cdot\tS=\emptyset$. By the Brouwer fixed point theorem we cannot have $\tB\subset\gamma\cdot\tB$ or  $\tB\supset\gamma\cdot\tB$, and since $\HH^3$ is non-compact we cannot have $\tB\cup\gamma\cdot\tB=\HH^3$. Hence 
$\tB\cap\gamma\cdot\tB=\emptyset$ for every $\gamma\in\Gamma-\{1\}$. It follows that $B:=p(\tB)$ is a $3$-ball with boundary $S$, and the proof of (\ref{it's irreducible}) is complete.

To prove (\ref{yes injective}), consider a torus $T\subset M$ which is incompressible.
The
image of the inclusion homomorphism $\pi_1(T)\to\pi_1(M)$ 
is a rank-$2$ free abelian subgroup of $\pi_1(M)$ which is defined up to conjugacy, and corresponds to a a rank-$2$ free abelian subgroup $X$ of $\Gamma$ which is also defined up to conjugacy. Since $\Gamma$ is discrete, $X$ must be parabolic. It is then a standard consequence of Shimizu's lemma (see for example \cite[Theorem 2.21]{series}) that there is a  horoball $H\subset\HH^3$, precisely invariant under $\Gamma$, whose stabilizer $\Gamma_H$ contains $X$. In particular $\Gamma_H$ is non-cyclic, and since it is parabolic and torsion-free it must be a rank-$2$ free abelian group. Hence $C_0:=H/\Gamma$ is a \cn. 
After possibly replacing $C_0$ by a smaller \cn, we may assume that $C_0\cap T=\emptyset$. Since  $X\le\Gamma_H$, and since $M$ is aspherical by (\ref{it's irreducible}), the inclusion map of $T$ into $M$ is homotopic in $M$  to a map whose image is contained in $C_0$. This inclusion map is therefore homotopic in $\overline{M-C_0}$  to a map whose image is contained in $\partial C_0$. It now follows from Lemma \ref{it's parallel}, applied with $N=\overline{M-C_0}$, that $T$ is boundary-parallel in $\overline{M-C_0}$;  that is, there is a submanifold $P$ of $M$, diffeomorphic to $T^2\times[0,1]$, such that $\partial P=T\cup\partial C_0$ and $P\cap C_0=\partial C_0$. Since $C_0$ is diffeomorphic to $T^2\times[0,\infty)$, the submanifold $C:=C_0\cup P$, which is closed as a subset of $M$ and has boundary $T$, is also diffeomorphic to $T^2\times[0,\infty)$. Furthermore, since $C_0$ has finite volume and $P$ is compact, $C$ also has finite volume.

To prove (\ref{not injective}), consider a torus $T\subset M$ which is not incompressible.
In this case, by \cite[Lemma 6.1]{hempel}, there is a disk $D\subset\inter M$ such that $D\cap T=\partial D$, and $\partial D$ is a non-trivial, and hence non-separating, curve on $T$. Let $E$ be a $3$-ball containing $D$, such that $A:=E\cap T$ is an annular neighorhood of $\partial D$ in $T$, and $A\subset\partial E$. Then $S:=(T\cup\partial E)-\inter A$ is a  $2$-sphere which must bound a  $3$-ball $B\subset M$, since $M$ is irreducible by (\ref{it's irreducible}). We must have either $E\subset B$ or $E\cap  B=\closure{(\partial E)-A}$. In the first case we have $T\subset B$, and (ii) holds. In the second case, $J:=E\cup B$ is a solid torus since $M$ is orientable, and $\partial J=T$, so that (b) holds.

\EndProof

%\Number\label{oldhockenapuss}
  %\redcomment{Give a discussion of Haken manifolds. Point out that if the boundary is %non-empty and the manifold is not a ball then by a result already quoted later on from %\cite{finiteness} we have $H_1(M;\QQ)\ne0$ and then $M$ contains an essential %(non-separating) surface by \cite[Lemma 6.6]{hempel}.
%}
%\EndNumber

The following well-known consequence of the sphere theorem is included, for example, in \cite[Theorem 8.2]{epstein}:

\Proposition\label{sphere consequence}An irreducible $3$-manifold with infinite fundamental group is aspherical.\NoProof\EndProposition

The following two facts about $3$-manifolds are well-known, but I have supplied proofs for completeness.

\Proposition\label{solid-torus}Let $N$ be a compact, irreducible, orientable
$3$-manifold such that every component of $\partial N$
is a torus. If $N$ is boundary-reducible then it is a solid torus.
\EndProposition

\Proof I will prove the corresponding statement in the PL category (see \ref{cat-a-gory}).
Let $N$ be a compact, irreducible, orientable PL
$3$-manifold such that every component of $\partial N$
is a torus. If $N$ is boundary-reducible, it contains an
essential disk $D$. By hypothesis the component of $\partial N$
containing $\partial D$ is a torus $T$. If $E$ denotes a regular
neighborhood of $D$ in $N$, then the component of
$\partial(\closure{N-E})$ that meets $T$ is a $2$-sphere. By
irreducibility  follows that $\closure{N-E}$ is a ball, so that
$N$ is a solid torus.
\EndProof

\Proposition\label{core}
Suppose that $M$ is an irreducible, orientable  $3$-manifold. Then there is a compact, irreducible submanifold $N$ of $M$ such that the inclusion homomorphism $\pi_1(N)\to\pi_1(M)$ is an isomorphism.
\EndProposition

\Proof
According to \cite{core}, there is a compact submanifold $N_0$ of $M$
such that the inclusion homomorphism $\pi_1(N_0)\to\pi_1(M)$ is an
isomorphism. By Lemma \ref{plug}, there is a compact
irreducible submanifold $N$ of $M$ such that $N\supset N_0$, and the
inclusion homomorphism $\pi_1(N_0)\to\pi_1(N)$ is surjective. It
follows that the inclusion homomorphism $\pi_1(N)\to\pi_1(M)$ is an
isomorphism. 
\EndProof

I will make use of the characteristic submanifold theory \cite{js}, \cite{johannson}. The information that I will need is summarized in the following statement:

\Proposition\label{absolutely}
Let $M$ be any Haken manifold. Then there is  there is a Seifert-fibered manifold $\Sigma\subset\inter M$, such that each component of $\partial\Sigma$ in incompressible in $M$, and having the following property: if $f:T^2\to M$ is any map such that $f_\sharp:\pi_1(T^2)\to\pi_1(M)$ is injective, then $f$ is homotopic to a map $g:T^2\to\inter M$ such that $g(T^2)\subset\Sigma$.
\EndProposition

\Proof
According to the statement of the Characteristic Pair Theorem on page 138 of \cite{js}, and the discussion of the case $T=\emptyset$ following that statement, there is a Seifert-fibered manifold $\Sigma\subset\inter M$, such that each component of $\partial\Sigma$ is incompressible in $M$, and having the following property: if $\cals$ is any Seifert fibered space, and $F:\cals\to M$ is any map which is nondegenerate (in the sense defined on p. 55 of \cite{js}, taking $\calf=\calt=\emptyset$), then $F$ is homotopic to a map whose image is contained in $\Sigma$.        Now if $f:T^2\to M$ is any map such that $f_\sharp:\pi_1(T^2)\to\pi_1(M)$ is injective, and if we let $q:T^2\times[0,1]\to T^2$ denote the projection to the first factor, it follows immediately from the definition that $f\circ q$ is a nondegenerate map of the Seifert fibered manifold $T^2 \times[0,1]$ into $M$. Hence $f\circ q$ is homotopic to a map of $T^2\times[0,1]$ into $M$ whose image is contained in $\Sigma$, and so $f$ is homotopic to a map $g:T^2\to\inter M$ such that $g(T^2)\subset\Sigma$.
\EndProof

(It may be shown that the submanifold $\Sigma$ given by Proposition \ref{absolutely} is unique up to ambient isotopy in $M$, but I will not need this fact. Note that although $\Sigma\subset\inter Q$,
every torus component $T$ of $\partial Q$ which is $\pi_1$-injective in $M$ is ``parallel'' to a component of $\partial\Sigma$. One may think of $\Sigma$ is an ``absolute'' characteristic submanifold of $Q$ which ``carries essential singular tori,'' as distinguished from the  ``relative'' characteristic submanifold which is defined only when $Q$ is boundary-irreducible, and carries both essential singular tori and essential singular annuli.)

The next two results, Propositions \ref{jaco} and \ref{life is never simple}, are closely related to Theorem VI.4.1 of \cite{js}. They are both essentially topological results, although I have found it more convenient to state Proposition \ref{life is never simple} in the setting of hyperbolic $3$-manifolds.

\Proposition\label{jaco}
Let $N$ be a compact, irreducible, orientable $3$-manifold such that $\pi_1(N)$ has rank $2$ and is not free. Then each component of $\partial N$ is a torus.
\EndProposition

\Proof
If $N$ is closed, there is nothing to prove. If $\partial N\ne\emptyset$ then $N$ has the homotopy type of a connected finite CW complex $K$ of dimension at most $2$. We may take $K$ to have only one $0$-cell. Let $m$ and $n$ denote, respectively, the numbers of $1$-cells and $2$-cells of $K$. Then $\pi_1(K)$ has a presentation with $m$ generators and $n$ relations. By definition, the {\it deficiency} of this presentation is $m-n$. It is shown in \cite{magnus} that if if $k$ is a positive integer, a finitely presented group that has rank $k$ and has a presentation of deficiency at least $k$ must be free. Since
$\pi_1(K)\cong\pi_1(N)$ has rank $2$ and is not free, the deficiency $m-n$ must be at most one. This gives
$$\frac12\chi(\partial N)=\chi(N)=\chi(K)=1-m+n\ge0.$$
Note that if some component of $\partial N$ were a sphere then by irreducibility $N$ would be a ball, which is impossible since $\pi_1(N)$ has rank $2$.
Thus $\partial N$ is a closed, orientable $2$-manifold with no sphere components and $\chi(N)\ge0$. Hence every component of $\partial N$ is a torus.
\EndProof

\Proposition\label{life is never simple}
Let $M$ be a hyperbolic $3$-manifold of infinite volume. Then every two-generator non-abelian subgroup of $\pi_1(M)$ is free.
\EndProposition

\Proof
Let $X\le \pi_1(M)$ be a two-generator non-abelian subgroup. Then $X\cong\pi_1(\tM)$ for some covering space $\tM$ of $M$. 
Since $M$ is irreducible by Assertion (\ref{it's irreducible}) of Proposition \ref{hypersimple}), it follows from \cite[p. 647, Theorem 3]{msy} that $\tM$ is irreducible.
By Proposition \ref{core}, there is a compact, irreducible submanifold $N$ of $M$ such that the inclusion homomorphism $\pi_1(N)\to\pi_1(M)$ is an isomorphism.
If $X\cong\pi_1(N)$ is not free, then it follows from Proposition \ref{jaco} that every component of $\partial M$ is a torus. 
Since $\pi_1(N)\cong X$ is non-abelian by hypothesis, $N$ is not a solid torus; hence by Proposition \ref{solid-torus}, $N$ is boundary-irreducible. If $T$ is any component of $\partial N$, it follows that $T$ is $\pi_1$-injective in $N$, and therefore in $M$ as well. Hence by Assertion (\ref{yes injective}) of Proposition \ref{hypersimple}, $T$ is the boundary of a submanifold $C_T$ of $M$ which is closed as a subset of $M$, is diffeomorphic to $T^2\times[0,\infty)$, and has finite volume. Since $\pi_1(N)$ is non-abelian, and since the inclusion homomorphism $\pi_1(N)\to\pi_1(M)$ is in particular injective, we cannot have $C_T\supset N$ for any $T$. Hence we must have $C_T\cap N=T$ for each $T$, so that $M=N\cup\bigcup_TC_T$. But this implies that $M$ has finite volume, a contradiction.
\EndProof

The following result has direct relevance to estimating Margulis numbers, and its corollaries were mentioned in the introduction.

\begin{proposition}\label{chained to the altar} Let $x$ and $y$ be non-commuting elements of 
$\isomplus(\HH^3)$ such that $\langle x,y\rangle$ is discrete and torsion-free and has infinite covolume. Then for every $P\in\HH^3$ we have
     $$\max(d(P,x\cdot P),d(P,y\cdot P))\ge\log3.$$   
\end{proposition}
 
\Proof
Set $\Gamma=\langle x,y\rangle$. Proposition \ref{life is never simple}, applied to the hyperbolic $3$-manifold $M=\HH^3/\Gamma$, shows that $\Gamma$ is free. The conclusion now follows from the case $k=2$ of \cite[Theorem 4.1]{surgery}.
\EndProof

\Corollary\label{born free}
Let $M$ be an orientable hyperbolic $3$-manifold such that every subgroup of rank at most $2$ in $\pi_1(M)$ has infinite index. Then $\mu(M)\ge\log3=1.09\ldots$.
\EndCorollary

\Proof
Write $M=\HH^3/\Gamma$, where $\Gamma\le\isomplus(\HH^3)$ is discrete and torsion-free. Let $x$ and $y$ be non-commuting elements of $\Gamma$. The hypothesis implies that $\langle x,y\rangle$ has infinite index in $\Gamma$, and hence has infinite covolume. By Proposition \ref{chained to the altar}, it follows that
$\max(d(P,x\cdot P),d(P,y\cdot P))\ge\log3$
for every $P\in\HH^3$. Hence $\log3$ is a Margulis number for $M$.
\EndProof

\Corollary\label{born to facebook}
Let $M$ be an orientable hyperbolic $3$-manifold such that either $H_1(M;\QQ)$ has rank at least $3$,   or 
$M$ is closed  and $H_1(M;\ZZ_p)$ has rank at least $4$ for some prime $p$. Then $\mu(M)\ge\log3$.
\EndCorollary

\Proof
If $H_1(M;\QQ)$ has rank at least $3$,   it is clear that every subgroup of rank at most $2$ in $\pi_1(M)$ has infinite index. If
$M$ is closed  and $H_1(M;\ZZ_p)$ has rank at least $4$ for some prime $p$, then by \cite[Proposition 1.1]{sw}, it is again true that every subgroup of rank at most $2$ in $\pi_1(M)$ has infinite index. Hence the result follows from Corollary \ref{born free}.
\EndProof

%  \redcomment{Use the following stuff as needed also:}

%{\bf A topological theorem}

%\begin{theorem}[Jaco-Shalen, Tucker]
%Let $M$ be a hyperbolic $3$-manifold (possibly with cusps and possibly of infinite %volume).   Let $J\le\pi_1(M)$ be a subgroup of rank at most two which has infinite %index in $\pi_1(M)$.   Then $J$ is either an abelian group or a free group of rank $2$.
%\end{theorem}

\section{An abstract bound for $\vol M$ when $\mu(M)<\log3$, and its consequences}
\label{abstract section}

Let $(\rho_n)_{n\ge1}$ be a sequence of representations of a group $X$ in  $\isomplus(\HH^3)$. Recall that the sequence $(\rho_n)$ is said to {\it converges algebraically} to a representation $\rho_\infty$ of $X$ in $\isomplus(\HH^3)$ if we have $\rho_n(\gamma)\to\rho_\infty(\gamma)$ for every $\gamma\in X$.

Recall that a subgroup of $\isomplus(\HH^3)$ is said to be {\it elementary} if it has an abelian subgroup of finite index. According to \cite[Proposition 2.1]{finiteness}, every torsion-free, elementary, discrete subgroup of $\isomplus(\HH^3)$ is abelian.

I will give an explicit proof of the following fact, which is implicit in \cite{jor}.

\Lemma\label{last one i hope} Let $(\rho_n)_{n\ge1}$ be a sequence of representations of a finitely generated group $\Phi$ in  $\isomplus(\HH^3)$. Suppose that $\rho_n(\Phi)$ is discrete, non-elementary  and torsion-free for each $n\in\NN$. Then there is a neighborhood $W$ of the identity in $\isomplus(\HH^3)$ such that $\tGamma_n\cap W=\{1\} $ for every $n\in\NN$.
\EndLemma

\Proof
If the conclusion is false, there is a sequence of elements $(x_n)_{n\ge1}$ of $\Phi$ such that $\rho_n(x_n)\ne1$ for $n=1,2,\dots$ but $\rho_n(x_n)\to\infty$ as $n\to\infty$. Since $\rho_n(\Phi)$ is a non-abelian, torsion-free, discrete subgroup of $\isomplus$, it has trivial center. In particular $\rho_n(x_n)\ne1$ is non-central in $\rho_n(\Phi)$ for each $n\in\NN$. Hence if $S$ is a finite generating set for $\Phi$,  then for each $n\in\NN$ there is an element $s_n\in S$ such that $\rho_n(x_n)$ and $\rho_n(s_n)$ do not commute. Since $S$ is finite, we may assume after passing to a subsequence that the $s_n$ are all the same element of $S$, say $s$. 

For each $n\in\NN$, the group $\Gamma_n:=\langle \rho_n(x_n),\rho_n(s)\rangle$ is contained in $\rho_n(\Phi)$, and is therefore discrete and torsion-free. Since $\rho_n(x_n)$ and $\rho_n(s)$ do not commute, $\Gamma_n$ is non-abelian, and is therefore non-elementary by \cite[Proposition 2.1]{finiteness}. 
If $X_n$ and $Y_n$ are elements of $\zzle(\CC)$ representing
$\rho_n(x_n)$ and $\rho_n(s)$ respectively, it then follows from
Jorgensen's inequality \cite[Lemma 1]{jor} that
\Equation\label{troels knows}
|(\trace X_n)^2-4|+|\trace(X_nY_nX_n^{-1}Y_n^{-1})-2|\ge1
\EndEquation
for each $n\in\NN$. As $n\to\infty$ we have 
$\rho_n(x_n)\to1$, while $\rho_n(s)\to\rho_\infty(s)$, where $\rho_\infty$ denotes the algebraic limit of the sequence $(\rho_n)$. After passing to a subsequence we may therefore assume that $X_n\to\pm I$ and that the sequence $(Y_n)$ has a limit. Hence $X_nY_nX_n^{-1}Y_n \to I$. It follows that the left hand side of (\ref{troels knows}) converges to $0$ as $n\ge\infty$, a contradiction.
\EndProof

The following result was stated in the introduction as Theorem A.

\begin{theorem}\label{abs tract}
Let $\lambda$ be a positive real number strictly less than $\log3$. Then there is a constant $V_\lambda$ such that for every orientable hyperbolic $3$-manifold $M$ with $\mu(M)<\lambda$ we have  $\vol M\le V_\lambda$ (and in particular $\vol M<\infty$).
\end{theorem}

\Proof
It follows immediately from Proposition \ref{chained to the altar} that if $M$ is any orientable $3$-manifold of infinite volume then $\log3$ is a Margulis number for $M$, so that
$\mu(M)\ge\log3>\lambda$.
 Hence if suffices to show that there is a constant $V_\lambda$ such that for every orientable hyperbolic $3$-manifold $M$ with $\infty>\vol M>V_\lambda$ we have $\mu(M)\ge\lambda$.

We reason by contradiction.   Assume there is a sequence $(M_n)_{n\ge1}$ of 
orientable finite-volume hyperbolic $3$-manifolds such that $\vol M_n\to\infty$ and  $\mu(M_n)<\lambda$, i.e. no $M_n$ admits $\lambda$ as a Margulis number.

For
each $n$ write $M_n=\HH^3/\Gamma_n$ for some torsion-free
cocompact discrete subgroup $\Gamma_n$ of $\isomplus(\HH^3)$.    
Then, by definition, for each $n$ there exist non-commuting elements
$x_n,y_n\in\Gamma_n$   and a point $P_n\in\HH^3$   such that 
$$\max(d(P_n,x_n\cdot P_n),d(P_n,y_n\cdot P_n))<\lambda.$$

After replacing each $\Gamma_n$ by a suitable conjugate of itself in
$\isomplus(\HH^3)$, we may assume that the $P_n$ are all the same point of $\HH^3$,    which I will denote by $P$.    Thus for each $n$ we have
\begin{equation}\label{maximillian}\max(d(P,x_n\cdot P),d(P,y_n\cdot P))<\lambda.\end{equation}

For each $n$, set $\tGamma_n:=\langle x_n,y_n\rangle$ and $\tM_n:=\HH^3/\tGamma_n$. Note that $\tGamma_n$ is discrete and torsion-free since $\Gamma_n$ is, and that $\tGamma_n$ is non-abelian---and hence non-elementary by \cite[Proposition 2.1]{finiteness}---since $x_n$ and $y_n$ do not commute.

Since $\lambda<\log3$, it follows from (\ref{maximillian}) and Proposition \ref{chained to the altar} that $\vol\tM_n<\infty$. On the other hand, $\tM_n$ covers $M_n$, and hence $\vol \tM_n\ge\vol M_n$.  In particular, $\vol\tM_n\to\infty$.

It follows from (\ref{maximillian}) that the $x_n$ and $y_n$ lie in a compact subset of $\isomplus(\HH^3)$. Hence, after passing to a subsequence, we may assume that the sequences $(x_n)$ and $(y_n)$ converge in $\isomplus(\HH^3)$ to limits $x_\infty$ and $y_\infty$. Again by (\ref{maximillian}), we have
\begin{equation}\label{welsh}
\max(d(P,x_\infty\cdot P),d(P,y_\infty\cdot P))\le\lambda.
\end{equation}

For $1\le n\le\infty$ we define a representation $\rho_n$ of the rank-$2$ free group $F_2=\langle\xi,\eta\rangle$ by $\rho_n(\xi)=x_n$, $\rho_n(\eta)=y_n$.  Thus $\rho_n(F_2)=\tGamma_n$ for each $n$. Since $\xi$ and $\eta$ generate $F_2$, and since $\rho_n(\xi)\to\rho_\infty(\xi)$ and $\rho_n(\eta)\to\rho_\infty(\eta)$ as $n\to\infty$, we have $\rho_n(\gamma)\to\rho_\infty(\gamma)$ for every $\gamma\in F_2$. By definition this means that the sequence $(\rho_n)$ converges algebraically to $\rho_\infty$.   Set $x_\infty=\rho_\infty(\xi)$, $y_\infty=\rho_\infty(\eta)$.

Let $D$ denote the set of representations of $F_2$ in $\isomplus(\HH^3)$ whose  images are discrete, torsion-free, and non-elementary.
According to \cite[Theorem 2.4]{finiteness} (a theorem essentially due to T. Jorgensen and P. Klein \cite{jk}), the limit of any algebraically convergent sequence of representations in $D$ is again in $D$.  Hence $\rho_\infty\in D$. Thus $\tGamma_\infty:=\rho_\infty(F_2)=\langle x_\infty,y_\infty\rangle$ is a discrete group.

According to \cite[Proposition 3.8]{jm}, since $(\rho_n)$ converges algebraically, the sequence of discrete groups $(\tGamma_n) $ has a geometrically convergent subsequence (in the sense defined in \cite{jm}). Hence without loss of generality we may assume that $(\tGamma_n) $ converges geometrically to some discrete group $\h\Gamma_\infty$. It then follows, again from \cite[Proposition 3.8]{jm},  that 
$\tGamma_\infty\le \h\Gamma_\infty$.

According to Lemma \ref{last one i hope}, there is a neighborhood $W$ of the identity in $\isomplus(\HH^3)$ such that $\tGamma_n\cap W=\{1\} $ for every $n\in\NN$. Let $E$ denote the set of all torsion-free subgroups $\Delta$ of $\isomplus(\HH^3)$ such that $\Delta\cap W=\{1\} $. (In particular each group in $E$ is discrete.) According to \cite[Theorem 1.3.1.4]{ceg}, $E$ is compact in the topology of geometric convergence. Since $\tGamma_n\in E$ for each $n\in\NN$, we have $\h\Gamma_\infty\in E$. In particular $\h\Gamma_\infty$ is torsion-free. We let $\h M_\infty$ denote the orientable hyperbolic $3$-manifold $\HH^3/\h\Gamma_\infty$.

Since $(\tGamma_n) $ converges geometrically to $\h\Gamma_\infty$, the sequence of orientable hyperbolic $3$-manifolds $(\tM_n)$ converges geometrically to $\h M_\infty$ in the sense of \cite[Chapter E]{bp}. If $\vol \h M_\infty$ were finite, it would then follow from \cite[Proposition E.2.5]{bp} that the sequence $(\vol\tM_n)$ had the finite limit $\vol \h M_\infty$, which contradicts 
$\vol\tM_n\to\infty$. Thus $\tM_\infty:=\HH^3/\rho_\infty(F_2)$ is a hyperbolic $3$-manifold of infinite volume. It therefore follows from Proposition \ref{chained to the altar} that
$$\max(d(P,x_\infty\cdot P),d(P,y_\infty\cdot P))\ge\log3.$$  
But this contradicts (\ref{welsh}).
\EndProof 
%Let $\Phi$ denote the subset of $D$ consisting of those discrete torsion-free representations whose images have finite covolume.   It is well known that the function $\rho\mapsto\vol(\HH^3/\rho(F_2))$ is continuous on $\Phi$.  If $\rho_\infty\in \Phi$, it follows that
%$$\vol(\HH^3/\rho_i(F_2))\to\vol(\HH^3/\rho_\infty(F_2)),$$
%  a contradiction since
%$$\vol(\HH^3/\rho_i(F_2))=\vol(\HH^3/\tGamma_i)\to\infty.$$  

\bigskip

The following two corollaries of Theorem \ref{abs tract} were also pointed out in the introduction.

\Corollary\label{synecdoche}
Let $\lambda$ be a positive real number strictly less than $\log3$. Then there is a 
there is a natural number $d_\lambda$ such that for every orientable hyperbolic $3$-manifold $M$ with $\mu(M)<\lambda$, the group $\pi_1(M)$ has a rank-$2$ subgroup of index at most $d_\lambda$.
\end{corollary}

\Proof
Let $v$ denote the infimum of the volumes of all
hyperbolic $3$-manifolds;  we have $v>0$, for example by
\cite[Theorem 1]{meyerhoff}. Let $V_\lambda$ be a positive real number having the property stated in Theorem \ref{abs tract}, and set 
$$d_\lambda=\lfloor \frac{V_\lambda}v\rfloor.$$

Suppose that $M$ is an orientable hyperbolic $3$-manifold 
such that $\mu(M)<\lambda$, i.e. such that $\lambda$ is not a Margulis number for $M$. Write $M_\infty=\HH^3/\Gamma$, where $\Gamma\le\isomplus(\HH^3)$ is discrete and torsion-free. Then by definition there exist a point $P\in\HH^3$ and non-commuting elements $x,y\in\Gamma$ such that
\Equation\label{Yes, AGAIN!}
\max(d(P,x\cdot P),d(P,y\cdot P))<\lambda.
\EndEquation 

Since $x$ and $y$ are non-commuting elements of $\tGamma$, it follows from (\ref{Yes, AGAIN!}) that $\lambda$ is not a Margulis number for $\tM$, i.e. that $\mu(M)<\lambda$. Hence $\infty\ge\vol \tM\le V_\lambda<\infty$, and since $\tM$ covers $M$ we have $\vol M<\infty$. Since $\vol M\ge v$, we find that
$$[\Gamma:\tGamma]=\frac{\vol\tM}{\vol M}\le
\frac{V_\lambda}v.$$
It follows that $[\Gamma:\tGamma]\le d_\lambda$, so that $\Gamma\cong\pi_1(M)$ has a rank-$2$ subgroup of index at most $d_\lambda$.
\EndProof

\Corollary\label{poor thing}
Let $\lambda$ be a positive real number strictly less than $\log3$. Then there is a 
natural number $k_\lambda$ such that for every orientable hyperbolic $3$-manifold $M$ with $\mu(M)<\lambda$, the group $\pi_1(M)$ has rank at most $k_\lambda$.
\end{corollary}

\Proof
Let $d_\lambda$ be a natural number having the property stated in Corollary \ref{synecdoche}. Set $k_\lambda=2+\lfloor\log_2d_\lambda\rfloor$. If $M$ is an orientable hyperbolic $3$-manifold with $\mu(M)<\lambda$, then the group $\Gamma=\pi_1(M)$ has a rank-$2$ subgroup $\tGamma$ of index at most $d_\lambda$. According to Proposition \ref{you peeked}, we have
$$\rank\Gamma\le \rank\tGamma+\log_2[\Gamma:\tGamma]\le2+d_\lambda$$
and hence $\pi_1(M)\cong \Gamma$ has rank at most $k_\lambda$.
\EndProof

\section{Margulis numbers and short relations}\label{short section}

\Notation\label{weird enough for ya?}
Let $\lambda$ be a real number $0<\lambda<(\log3)/2$. Then for any sufficiently large positive integer $N$ we have
\Equation\label{you're so weird}
\frac{3^{N}-1}{4N+1}\ge2667(\sinh(2N\lambda+.104)-(2N\lambda+.104)).
\EndEquation
(The natural logarithm of the left hand side of (\ref{you're so weird}) is asymptotic to $N\log3$, whereas the natural logarithm of the right hand side is asymptotic to $2N\lambda>N\log3$.)

I shall denote by $N(\lambda)$ denote the smallest positive integer $N$ for which (\ref{you're so weird}) holds.
\EndNotation
 
Here is one simple estimate of $N(\lambda)$:

\Proposition\label{nestimate}
Let $\lambda\in(0.1,(\log3)/2)$ be given, and set $\beta=1/((\log3)-2\lambda)$. Then
\Equation\label{real nestimate}
N(\lambda)<1+110\beta\log\beta.
\EndEquation
\EndProposition

\Proof
Since $\lambda>0.1$ we have $\beta>1/((\log3)-0.2)>1.11$. Hence the right hand side of (\ref{real nestimate}) is bounded below by $1+(110)(1.11)(\log 1.11)=13.7\ldots$. We may therefore assume without loss of generality that $N(\lambda)\ge14$.

Set $n=N(\lambda)-1\ge13$. From the definition of $N(\lambda)$ we have
\Equation\label{obversely weird}
\frac{3^{n}-1}{4n+1}<2667(\sinh(2n\lambda+.104)-(2n\lambda+.104)).
\EndEquation
Since in particular we have $n\ge3$, the left hand side of (\ref{obversely weird}) is bounded below by $3^n/5n$. The right hand side is obviously bounded above by $2667\exp(2n\lambda+.104)/2$. Hence
$$3^n\le\frac{2667}2\cdot5n\cdot\exp(2n\lambda+.104)<7400ne^{2n\lambda},$$
which upon taking logarithms and using the definition of $\beta$ gives
\Equation\label{poultroon}
n<\beta\log(7400n).
\EndEquation
Now suppose that (\ref{real nestimate}) is false, so that $n\ge110\beta\log\beta$. Define a function $g(x)$ for $x>0$ by $g(x)=x/\log(7400x)$. Then $g(x)$ is monotone increasing for $x>e/7400$, and since $n\ge110\beta\log\beta>12.7$, we have 
$g(n)\ge g(110\beta\log\beta)$. With (\ref{poultroon}) this gives
$$\frac{110\beta\log\beta}{\log(7400\cdot 110\beta\log\beta)} \le \frac{n}{\log(7400n)}<\beta,$$
so that
$${110\log\beta}<\log7400+\log110+\log\beta+\log\log\beta<13.61+\log\beta+\log\log\beta,$$
i.e.
\Equation\label{hermione}
{109\log\beta}-\log\log\beta<13.61.
\EndEquation
On the other hand, the function $h(x):=109x-\log x$ is monotonically increasing for $x>1/109$. Since $\log\beta>\log 1.11=0.104\ldots>1/109$, we have
$${109\log\beta}-\log\log\beta=h(\log\beta)>h(\log1.11)=13.63\ldots,$$
which contradicts (\ref{hermione}).
\EndProof

\begin{proposition}\label{first one}
Let
$M$ be an orientable hyperbolic $3$-manifold, and write $M=\HH^3/\Gamma$, where $\Gamma\le\isomplus(\HH^3)$ is discrete and torsion-free. Let $\lambda<(\log3)/2$ be given, and let $x$ and $y$ be non-commuting elements of $\Gamma$ such
that $\max(d(P,x\cdot P),d(P,y\cdot P))<\lambda$.    Then there is a reduced word $W$ in two letters, with $0<\length  W\le 8N(\lambda)$, such that $W(x,y)=1$. (Here $N(\lambda)$ is defined by  \ref{weird enough for ya?}.)  
\end{proposition}

\Proof
Set $\mu=0.104$. I pointed out in the introduction that according to \cite{meyerhoff}, we have $\mu<\mu_+(3)$; that is, $\mu$ is a Margulis number for every orientable hyperbolic $3$-manifold.  Set $N=N(\lambda)$. 

Since $\mu$ is in particular a Margulis number for $M$, the elements $\gamma\in\Gamma$ such that $d(\gamma\cdot P,P)<\mu$ generate an abelian subgroup $C$ of $\Gamma$. 

Let $F_2$ denote a free group on two (abstract) generators $\xi$ and $\eta$. We identify $F_2$ with the set of all reduced words in $\xi$ and $\eta$, so that $V(\xi,\eta)=V$ for every reduced word $V$. Let $\phi:F_2\to\Gamma$ denote the unique homomorphism such that $\phi(\xi)=x$ and $\phi(\eta)=y$; then $\phi(V)=V(x,y)$ for every reduced word $V$.

For every positive integer $n$, let $\calv_n\subset F_2$ denote the set of all reduced words of length at most $n$ in $\xi$ and $\eta$.
If $m$ and $n$ are positive integers, then for all $V\in\calv_m$ and $V'\in\calv_n$, we can concatenate $V$ and $V'$ and then reduce the resulting word to obtain a reduced word of length at most $m+n$ representing the product $VV'$ in $F_2$. Hence
\Equation\label{give up the chainu8}
\calv_m\calv_n\subset\calv_{m+n}\text{ for all }m,n>0.
\EndEquation  
Note also that
\Equation\label{hu 8 the chain}
\calv_n^{-1}=\calv_n\text{ for every }n>0.
\EndEquation

For any $k>0$, the number of reduced words in $\xi$ and $\eta$ of length exactly $k$ is $4\cdot3^{k-1}$. Summing from $k=1$ to $k=n$, we deduce that
\Equation\label{sum like it dim}
\#(\calv_n)-1=2(3^{n}-1) \text{ for every }n\ge1.
\EndEquation
I will assume that: 
\Equation\label{o'galavant}
\phi(W)\ne1\text{ for every }W\in\calv_{ 8N}-\{1\}.
\EndEquation
Under the assumption (\ref{o'galavant}) I will derive a contradiction, thereby proving the proposition.

Assuming (\ref{o'galavant}), I claim:
\Equation\label{pogo shtick}
\#(\calv_N\cap\phi^{-1}(\gamma C))\le4N+1 \text{ for any left   coset }\gamma C\text{ of }C\text{ in }\Gamma.
\EndEquation

To prove (\ref{pogo shtick}) I will first consider the special case in which $\phi^{-1}(C)\cap\calv_{2N}=\{1\}$. In this case, if $V$ and $V'$ are elements of $\calv_N\cap\phi^{-1}(\gamma C)$ for a given left   coset $\gamma C$ of $C$ in $\Gamma$, we have $U:=V^{-1}V'\in\phi^{-1}(C)$; but by (\ref{give up the chainu8}) and (\ref{hu 8 the chain}) we have $U\in\calv_{2N}$. The assumption  
$\phi^{-1}(C)\cap\calv_{2N}=\{1\}$ then implies that $U=1$ and hence that $V=V'$. Hence in this case we have $\#(\calv_N\cap\phi^{-1}(\gamma C))\le1$, which is stronger than (\ref{pogo shtick}).

Hence in proving  (\ref{pogo shtick}) we may assume that 
$\phi^{-1}(C)\cap\calv_{2N}\not=\{1\}$.
Of course we may also assume that $\calv_N\cap\phi^{-1}(\gamma C)\not=\emptyset$. Let us fix an element $U_0\ne1$ of $\phi^{-1}(C) \cap\calv_{2N}$ and an element $V_1$ of $\calv_N\cap\phi^{-1}(\gamma C)$. Let $\hC$ denote the centralizer of $U_0$ in  $F_2$; then $\hC$ is cyclic, since $F_2$ is free and $U_0\ne1$.

Let us define an injective map $J:\calv_N\cap\phi^{-1}(\gamma C)\to F_2$ by $J(V)=V^{-1}V_1$. Since
$V_1\in\phi^{-1}(\gamma C)$, we have $J(V)\in\phi^{-1} (C)$ for every $V\in\calv_N\cap\phi^{-1}(\gamma C)$. Since $V_1\in\calv_N$, it follows from  (\ref{give up the chainu8}) and (\ref{hu 8 the chain}) that $J(V)\in\calv_{2N}$ for every $V\in\calv_N\cap\phi^{-1}(\gamma C)$. Thus $J$ maps $\calv_N\cap\phi^{-1}(\gamma C)$ into $\calv_{2N}\cap\phi^{-1}(C)$. 

Let $V\in\calv_N\cap\phi^{-1}(\gamma C)$ be given.
Since $U_0,J(V)\in\phi^{-1}(C)$, we have $\phi(J(V)U_0J(V)^{-1}U_0^{-1}=1$. Since $J(V)$ and $U_0$ belong to $\calv_{2N}$, it follows from (\ref{give up the chainu8}) and (\ref{hu 8 the chain}) that 
$J(V)U_0J(V)^{-1}U_0^{-1}\in\calv_{8N}$. Applying (\ref{o'galavant}) with $W=J(V)U_0J(V)^{-1}U_0^{-1})$ we deduce that $W=1$, i.e. that $J(V)$ commutes with $U_0$. This means that $J(V)\in\hat C$. 

Thus $J$ is an injection from $\calv_N\cap\phi^{-1}(\gamma C)$ to $\calv_{2N}\cap\hC$, and so $\#(\calv_N\cap\phi^{-1}(\gamma C))\le \#(\calv_{2N}\cap\hC)$. But according to Lemma \ref{o hula who}, appplied with $k=2N$ and with $Z=\hC$, we have $\#(\calv_{2N}\cap\hC)\le4N+1$. Thus (\ref{pogo shtick}) is established.

Now let $\call$ denote the set of all left cosets of $C$ in $\Gamma$, and define a map $\psi:F_2\to\call$ by $\psi(V)=\phi(V)C$. We may then paraphrase (\ref{pogo shtick}) by saying that the fibers of the surjection $\psi|\calv_N:\calv_N\to\psi(\calv_N)$ have cardinality at most $4N+1$. If we set $r=\#(\psi(\calv_N))$, it follows that $r\ge\#(\calv_N)/(4N+1)$. Combining this with (\ref{sum like it dim}), we find that
\Equation\label{barnacles to you}
r>\frac{2(3^{N}-1)}{4N+1}.
\EndEquation

Let $b$ denote the open ball of radius $\mu/2$ centered at $P$, and let 
$B$ denote the ball of radius $N\lambda+(\mu/2)$ centered at $P$. We have
\Equation\label{volare o o}
\vol b=\pi(\sinh(\mu)-\mu)= 0.000589\ldots
\EndEquation
and
\Equation\label{o o o o}
\vol B=\pi(\sinh(2N\lambda+\mu)-(2N\lambda+\mu)).
\EndEquation

Since 
$d(x\cdot P, P)<\lambda$ and $d(y\cdot P, P)<\lambda$, and since $x$ and $y$ are isometries, we have 
$d(\phi(V)\cdot P, P) =d(V(x,y)\cdot P, P)<N\lambda$
for every $V\in\calv_N$. It follows that
\Equation\label{if i had them i'd be king}
\phi(V)\cdot b\subset B\text{  for every }V\in\calv_N.
\EndEquation
According to the definition of $r$, there are elements $\gamma_1,\ldots,\gamma_r$ of $\phi(\calv_N)$ which represent distinct left cosets of $C$ in $\gamma$. From (\ref{if i had them i'd be king}) we have
\Equation\label{faith and begorrah}
\gamma_i\cdot b\subset B\text{  for }i=1,\ldots,r.
\EndEquation
If $i$ and $j$ are distinct indices in $\{1,\ldots,r\}$ we have $\gamma_j^{-1}\gamma_i\notin C$, which by the definition of $C$ gives $
d(\gamma_i\cdot P,\gamma_j\cdot P,P)=
d(\gamma_j^{-1}\gamma_i\cdot P,P)\ge\mu$, so that
\Equation\label{whatcha doin' in disjoint}
\gamma_i\cdot b\cap \gamma_j\cdot b=\emptyset\text{ for all distinct indices }i\text{ and }j\text{ in }\{1,\ldots,r\}.
\EndEquation
From (\ref{faith and begorrah}) and (\ref{whatcha doin' in disjoint}) it follows that $r\vol b\le\vol B$. Combining this with (\ref{barnacles to you}), (\ref{volare o o}) and (\ref{o o o o}), we obtain
$$\begin{aligned}
\frac{2(3^{N}-1)}{4N+1}&< r\le\frac{\vol B}{\vol b}\cr&\le {\pi(\sinh(2N\lambda+\mu)-(2N\lambda+\mu))}/{0.000589}\cr
&<5334
(\sinh(2N\lambda+\mu)-(2N\lambda+\mu))
\end{aligned}
$$
which is a contradiction, since by definition (\ref{you're so weird}) holds with $N=N(\lambda)$. 
\EndProof

\section{From a short relation to a volume bound}\label{short-bounded section}

The main result of this section is Proposition \ref{short to bounded}. As I mentioned in the introduction, the basic method of proof of that proposition is due to Cooper \cite{cooper}. Among the several preliminaries results needed to apply Cooper's method in the present situation, Lemma \ref{oh far out} is the deepest, while Proposition \ref{wow what a variety of characters} seems to be of particular independent interest.

\begin{lemma}\label{same old characters}Let $N$ and $Q$ be irreducible, orientable $3$-manifolds. Suppose that $\pi_1(N)$ is isomorphic to a subgroup of $\pi_1(Q)$ and that $ N$ is closed. Then either $Q$ is closed, or $N$ is simply connected. \end{lemma}
 
\Proof
First consider the case in which $\pi_1(N)$ is infinite. In this case $\pi_1(Q)$ is also infinite, and $Q$ and $N$ are aspherical by Proposition \ref{sphere consequence}.  Fix a base point $q$ in $Q$ and a subgroup $J$ of $\pi_1(Q,q)$ isomorphic to $\pi_1(N)$, and let $(\tQ,\tq)$ denote the based covering of $(Q,q)$ determined by $J$. Since $N$ and $\tQ$ are aspherical and have isomorphic fundamental groups, they are homotopy-equivalent. If $ Q$ is not closed then $ \tQ$ is not closed, and hence $H_3(\tQ;\ZZ)=0$. But $H_3(N;\ZZ)\ne0$ since $N$ is closed, and we have a contradiction to the homotopy-equivalence of $\tQ$ and $N$. Hence in this case $ Q$ must be closed.

Now consider the case in which $\pi_1(N)$ is finite. In this case I will assume that $\pi_1(N)$ is non-trivial and show that $Q$ is closed, thus establishing the conclusion. The assumption that $\pi_1(N)$ is finite and non-trivial implies that $\pi_1(Q)$ has torsion, and so $Q$ cannot be aspherical \cite[Lemma 8.4]{epstein}. By 
Proposition \ref{sphere consequence}, $\pi_1(Q)$ is finite.  
We apply Proposition \ref{core}, letting $Q$ play the role of $M$ in that proposition. This gives a compact, irreducible submanifold $Q_0$ of $Q$ such that the inclusion homomorphism $\pi_1(Q_0)\to\pi_1(Q)$ is an isomorphism.
In particular, $\pi_1(Q_0)$ is finite, so that $H_1(Q_0;\QQ)=0$, and hence every component of $\partial Q_0$ is a sphere (for example by \cite[Proposition 2.2]{finiteness}).
Since $Q_0$ is irreducible it follows that either $\partial Q_0=\emptyset$ or $Q_0$ is a ball. In the latter case the hypothesis gives $\pi_1(N)=\{1\}$, a contradiction. Hence $\partial Q_0=\emptyset$, which implies that $Q=Q_0$ and that $Q$ is closed.
\EndProof

\begin{proposition}\label{wow what a variety of characters}Let $N$ and $Q$ be irreducible, orientable $3$-manifolds (possibly with boundary). Suppose that $N$ is compact, that $\pi_1(N)$ is non-cyclic and that every component of $\partial N$ is a torus. Let $h:N\to Q$ be a map such that $h_\sharp:\pi_1(N)\to\pi_1(Q)$ is injective. Then there exist a compact submanifold $Q_0$ of $Q$ such that every component of $\partial Q_0$ is a torus, and a map $h_0:N\to Q$ homotopic to $h$, such that $h_0(N)\subset Q_0$.
\end{proposition}

\Proof
First note that since $N$ is irreducible, $\pi_1(N)$ is non-cyclic, and every component of $\partial N$ is a torus, it follows from Proposition \ref{solid-torus} that $N$ is boundary-irreducible. 

First consider the case in which $N$ is closed. Note that since $h_\sharp$ is injective,
$\pi_1(N)$ is isomorphic to a subgroup of $\pi_1(Q)$. Since $\pi_1(N)$ is non-cyclic, it then follows from Lemma \ref{same old characters} that $Q$ is closed. Hence the conclusion holds in this case if we set $Q_0=Q$.

Next consider the case in which $\partial N\ne\emptyset$ and $Q$ is compact. Let us fix a Seifert-fibered space $\Sigma\subset\inter M$ having the property stated in Proposition \ref{absolutely}. Let us denote by $J$ the union of all components $C$ of $\overline{Q-\Sigma}$ such that $\partial C\subset\partial \Sigma$.

In this case I will set $Q_0=J\cup\Sigma$, so that every component of $\partial Q_0$ is a torus, and I will construct a map $h_0:N\to Q$, homotopic to $h$, such that $h_0(N)\subset Q_0$.

Consider an arbitrary component $T$ of $\partial N$. By hypothesis $T$ is a torus.  Since $N$ is boundary-irreducible, $T$ is $\pi_1$-injective in $N$. Since $h_\sharp:\pi_1(N)\to\pi_1(Q)$ is injective it follows that $(h|T)_\sharp:\pi_1(T)\to\pi_1(Q)$ is injective. It then follows from the property of $\Sigma$ stated in Proposition \ref{absolutely} that $h|T:T\to Q$ is homotopic to a map whose image is contained in $\Sigma$. Since this is true for each component $T$ of $\partial N$, it follows that $h|\partial N$ is homotopic to a map $g:\partial N\to Q$ such that $g(\partial N)\subset\inter\Sigma$. By the homotopy extension property of polyhedra, $g$ extends to a map from $N$ to $Q$ which is homotopic to $h$.

It follows from \cite[Lemma 6.5]{hempel} that we can choose an extension $h_1:N\to Q$ of $g$, homotopic to $h$, so that $h_1$ is transverse to $\partial \Sigma$ and so that each component of $h_1^{-1}(\partial\Sigma)$ is incompressible. Since $h_1(\partial N)=g(\partial N)\subset\inter\Sigma$, we have $h_1^{-1}(\partial\Sigma)\subset\inter N$.

I claim:
\Claim\label{talmudic} If $K$ is a component of $h_1^{-1}(\overline{Q-Q_0})$, then $h_1|K$ is homotopic rel $\partial K$ to a map $h_K$ with $h_K(K)\subset\partial\Sigma$. (In particular
$h_1(\partial K)\subset\partial\Sigma$.)
\EndClaim
To prove \ref{talmudic}, consider a component $K$  of $h_1^{-1}(\overline{Q-Q_0})$, and let $C$ denote the component   of $\overline{Q-Q_0}$ containing $h_1(K)$. 
Since $h_1$ is transverse to $\partial\Sigma$ and maps $\partial N$ into $\inter\Sigma$, we have $h_1(\partial K)\subset C\cap\partial\Sigma$. In particular $h_1|K:K\to C$  is a boundary-preserving map.

Every component of the frontier of $K$ in $N$ is a component of $h_1^{-1}(\partial\Sigma)$ and is therefore incompressible in $N$. Hence $K$ is $\pi_1$-injective in $N$. 
On the other hand, $h_\sharp:\pi_1(N)\to\pi_1(Q)$ is injective by hypothesis, and since $h_1$ and $h$ are homotopic, $(h_1)_\sharp:\pi_1(N)\to\pi_1(Q)$ is also injective. This shows that $(h_1|K)_\sharp:\pi_1(K)\to\pi_1(Q)$ is injective, and so in particular
$(h_1|K)_\sharp:\pi_1(K)\to\pi_1(C)$ is injective.

Since $N$ is irreducible and boundary-irreducible, and since we have observed that the components of the frontier of $K$ in $N$ are incompressible, the manifold $K$ is also irreducible and boundary-irreducible. Since the components of $\partial N$ are tori and the components of the frontier of $ K$ are incompressible, $\partial K$ has no sphere components.
Since $\partial N\ne\emptyset$, we have $\partial K\ne\emptyset$. According to \cite[Lemma 6.7]{hempel}, it follows that $K$ is a Haken manifold. Since $K$ is boundary-irreducible it is not a solid torus. It now follows from 
\cite[Theorem 13.6]{hempel} that every boundary-preserving map from $K$ to $C$ which induces an injection of fundamental groups  is homotopic rel $\partial K$ either to a covering map or to a map whose image is contained in $\partial C$.

Suppose that
$h_1|K:K\to C$ is homotopic rel $\partial K$ to a covering map. Then in particular we have $h(\partial K)=\partial C$. Since we have seen that $h_1(\partial K)\subset C\cap\partial\Sigma$, this implies that $\partial C\subset\partial\Sigma$, which by the definition of $Q_0$ implies that $C\subset Q_0$. This is a contradiction, since
$C$ is a component   of $\overline{Q-Q_0}$. Hence $h_1|K:K\to C$ is homotopic rel $\partial K$ to a map whose image is contained in $\partial C$. Since $h_1(\partial K)\subset C\cap\partial\Sigma$, it follows that we in fact have $h_K(\partial K)\subset C\cap\partial\Sigma$. This proves \ref{talmudic}.

We may now define a map $h_0:N\to Q$ by letting $h_0$ agree with $h_1$ on $h_1^{-1}(Q_0)$, and setting $h_0|K=h_K$ for each component $K$ of $h_1^{-1}(\overline{Q-Q_0})$. It is immediate from the properties of the $h_K$ stated in \ref{talmudic} that $h_0$ is well-defined, is homotopic to $h_1$ (and hence to $h$), and maps $N$ into $Q_0$. This establishes the conclusion in this case.

There remains the case in which $\partial N\ne\emptyset$ and $Q$ is non-compact. In this case, fix a compact submanifold $R_0$ of $Q$ containing $h(N)$. By Lemma \ref{plug}, there is a compact
irreducible submanifold $R$ of $M$ such that $R\supset R_0$. The
hypotheses of the propositon continue to hold if $Q$ is replaced by
$R$. As $R$ is compact, the case of the proposition already proved
gives  a compact submanifold $Q_0$ of $R\subset Q$ such that every
component of $\partial Q_0$ is a torus, and a map $h_0:N\to R$
homotopic to $h$ in $R$ (and hence in $Q$), such that $h_0(N)\subset
Q_0$. \EndProof

%  \redcomment{Be careful about what category you're in, e.g. in the following statement.}

Although it is convenient to state the following result in the hyperbolic setting, the proof is essentially topological.

\begin{lemma}\label{oh far out}Let $M$ be a hyperbolic $3$-manifold, let $K$ be a compact, path-connected space such that $\pi_1(K)$ has rank $2$ and is not free, and let $f:K\to M$ be a continuous map such that $f_\sharp:\pi_1(K)\to\pi_1(M)$ is surjective and $f(K)\subset M$ is a polyhedron with respect to some $C^1$ triangulation of $M$.   Then for each component $C$ of $M-f(K)$, the  image of the inclusion homomorphism $\pi_1(C)\to\pi_1(M)$ is abelian.
\end{lemma}

\Proof 
We may assume that $\pi_1(M)$ is non-abelian, as otherwise the conclusion is obvious.

We fix a PL structure on $M$, defined by a $C^1$ triangulation, in which $f(K)$ is a polyhedron.

For the purpose of this proof, a compact polyhedron in $M$ will be termed ``small'' if it is contained in a  $3$-ball in $M$. Let $R$ be a regular neighborhood of $f(K)$ in $M$. Let $Q$ denote the union of $R$ with all small components of $\overline{M-R}$. Then no component of 
$\overline{M-Q}$ is small. If $S$ is any  $2$-sphere in $\inter Q$, then $S$ bounds a  $3$-ball $B\subset M$, since $M$ is irreducible by Assertion (\ref{it's irreducible}) of Proposition \ref{hypersimple}. If $B\not\subset Q$, then $B$ contains a component of 
$\overline{M-Q}$, which by definition must be small, a contradiction; hence $B\subset Q$. This shows that $Q$ is irreducible.

Fix a base point $x\in K$, and set $q=f(x)\in f(K)\subset R\subset Q$. We may regard $f$ as a map from the based space $(K,x)$ to the based space $(Q,q)$. Let $I\le\pi_1(Q,q)$ denote the image of the homomorphism $f_\sharp:\pi_1(K,x)\to\pi_1(Q,q)$. Then $I$ determines a based covering space $p:(\tQ,\tq)\to(Q,q)$, and there is a unique lift $\tf:K\to\tQ$ such that $\tf(x)=\tq$. 
If $i:Q\to M$ denotes inclusion, we have a commutative diagram
$$
\begin{xy}
(0,0)*+{\pi_1(K,x)}="a";(30,20)*+{\pi_1(\tQ,\tq)}="b";(30,0)*+{\pi_1(Q,q)}="c";
(60,0)*+++{\pi_1(M,q).}="d";
{\ar^{\tf_\sharp}"a";"b"}
{\ar^{f_\sharp}"a";"c"}
{\ar^{p_\sharp}"b";"c"}
{\ar^{i_\sharp}"c";"d"}
\end{xy}
$$
According to the hypothesis, $f_\sharp\circ i_\sharp$ is surjective, and hence by commutativity of the diagram, $i_\sharp\circ p_\sharp:\pi_1(\tQ,\tq)\to \pi_1(M,q)$ is surjective. Since $\pi_1(M)$ is non-abelian, it follows that $\pi_1(\tQ)$ is non-abelian, and in particular non-cyclic. On the other hand, the construction of $(\tQ,\tq)$ implies that 
$\tf_\sharp$ is surjective. Thus the non-cyclic group $\pi_1(\tQ)$ is a homomorphic image of $\pi_1(K)$, which by hypothesis is a non-free group of rank $2$. It follows that $\pi_1(\tQ)$ is also a non-free group of rank exactly $2$.

Since $Q$ is irreducible, it follows from \cite[p. 647, Theorem 3]{msy} that $\tQ$ is irreducible.
We now apply Proposition \ref{core}, letting $\tQ$ play the role of $M$ in that proposition. This gives a compact, irreducible submanifold $N$ of $\tQ$ such that the inclusion homomorphism $\pi_1(N)\to\pi_1(\tQ)$ is an isomorphism. In particular, $\pi_1(N)$ is a non-free group of rank exactly $2$. Hence by Proposition \ref{jaco}, each component of $\partial N$ is a torus.

Set $h=p|N:N\to Q$. Since $p_\sharp:\pi_1(\tQ)\to\pi_1(Q)$ is injective, and since the inclusion homomorphism $\pi_1(N)\to\pi_1(\tQ)$ is an isomorphism. We have seen that
$N$ and $Q$ are compact, irreducible, orientable $3$-manifolds that $\pi_1(N)$ has rank $2$, and is in particular non-cyclic; and that every component of $\partial N$ is a torus. Thus all the hypotheses of Proposition \ref{wow what a variety of characters} hold. The latter proposition now gives a compact submanifold $Q_0$ of $Q$ such that every component of $\partial Q_0$ is a torus, and a map $h_0:N\to Q$ homotopic to $h$, such that $h_0(N)\subset Q_0$.

We have seen that $i_\sharp\circ p_\sharp:\pi_1(\tQ)\to\pi_1(M)$ is surjective. Since the inclusion homomorphism $\pi_1(N)\to\pi_1(\tQ)$ is an isomorphism, it follows that
$i_\sharp\circ h_\sharp=i_\sharp\circ (p|N)_\sharp :\pi_1(N)\to\pi_1(M)$ is surjective. Since the maps $h_0,h:N\to Q$ are homotopic, $i_\sharp\circ (h_0)_\sharp:\pi_1(N)\to\pi_1(M)$ is also surjective. Since $h_0(N)\subset Q_0$, the inclusion homomorphism $\pi_1(Q_0)\to M$ is surjective. 

To establish the conclusion of the lemma, it suffices to show that for every component $c$ of $\closure{M-R}$, the image of the inclusion homomorphism $\pi_1(c)\to\pi_1(M)$ is abelian. According to the definition of $Q$, any such $c$ is either a small component of $\overline{M-R}$, in which case the image of the inclusion homomorphism is trivial, or a component of $\overline{M-Q}$. In the latter case, $c$ is contained in a component $c_0$ of $\overline{M-Q_0}$, and I shall complete the proof by showing that the image of the inclusion homomorphism $\pi_1(c_0)\to\pi_1(M)$ is abelian. Let us fix any component $T$ of $\partial c_0$. Then $T$ is a component of $\partial Q_0$ and is therefore a torus. It follows from Assertions (\ref{yes injective}) and (\ref{not injective}) of Proposition \ref{hypersimple} that $T$ is the boundary of a $3$-dimensional submanifold $J$ of $M$, closed as a subset of $M$, such that the inclusion homomorphism $\pi_1(J)\to\pi_1(M)$ has abelian image. Since $Q_0$ is connected, we must have either $J=c_0$, in which case the required conclusion holds, or $J\supset Q_0$. The latter alternative would imply that the inclusion homomorphism $\pi_1(Q_0)\to\pi_1(M)$ has abelian image. This is impossible, as this inclusion homomorphism is surjective and $\pi_1(M)$ is non-abelian.
\EndProof

\begin{lemma}\label{eh golly you}Let $M$ be a hyperbolic $3$-manifold of finite volume, and let $C$ be a compact $3$-dimensional (smooth) submanifold of $M$.   Suppose that the inclusion homomorphism $\pi_1(C)\to\pi_1(M)$ has abelian image.   
Then
$$\vol C\le\frac12\area \partial C.$$
\end{lemma}

\Proof 
This is essentially Lemma 4.4 of \cite{A-L}. The authors of \cite{A-L} assume $C$ to be ``PL'' in a sense that it is not clear to me, but their proof can be read in the smooth category without change.
\EndProof

\begin{lemma}\label{jeweler you failed}Let $M$ be a hyperbolic $3$-manifold of finite volume. Let $X\subset M$ be a compact set which is a $2$-dimensional polyhedron with respect to some smooth triangulation of $M$.   Suppose that for every component $C$ of $M-X$, the inclusion homomorphism $\pi_1(C)\to\pi_1(M)$ has abelian image.   
Then
$$\vol M\le\area X.$$
\end{lemma}

\Proof
Let $m$ denote the number of cusps of $M$.
Let $\epsilon>0$ be given. 
Since $M$ has finite volume, there are disjoint \cn s $V_1,\ldots,V_m$ in $M$, for some $m>0$, such that $M_0:=\overline{M-(V_1\cup\cdots V_m)}$ is compact. Set $V=V_1\cup\cdots\cup V_m$. After replacing the $V_i$ by smaller \cn s if necessary, we may assume that each $V_i$ has volume less than $\epsilon/m$, that
each $\partial V_i$ has area less than $\epsilon/m$, and that the $V_i$ are disjoint from $X$. Hence $\vol V<\epsilon$, $\area\partial V<\epsilon$, and $X\subset\inter M_0$.

Since $X$ is a $2$-dimensional polyhedron with respect to the given smooth triangulation of $M$, there is an open neighborhood $U$ of $X$ in $\inter M_0$ with $\vol U<\epsilon$. Let $N\subset U$ be a smooth regular neighborhood of $X$ in the sense of \cite{hirsch}. We may choose $N$ in such a way that $\area\partial N<\epsilon+2\area X$.
% \redcomment{Justification needed?}

Let $C_1,\ldots,C_n$ denote the components of $\overline{M_0-N}$. Then each $C_j$ is contained in a component of $M-X$, and hence the inclusion homomorphism $\pi_1(C_j)\to\pi_1(M)$ has abelian image. 
According to Lemma \ref{eh golly you}, we have
$$
\vol C_j\le\frac12\area \partial C_j
$$
for $j=1,\ldots,n$. Hence
$$\begin{aligned}
\vol(\overline{M_0-N})&\le\frac12\area \partial(\overline{M_0-N})\cr
&=\frac12\area\partial M_0+\frac12\area\partial N.
\end{aligned}
$$
Since $\area\partial M_0=\area\partial V<\epsilon$ and $\area\partial N<\epsilon+2\area X$, it follows that
$$\vol(\overline{M_0-N})<\epsilon+\area X.$$
%But since $X\subset X\cup\partial V$, we have
%$$\area X\le\area X+\area\partial V<\epsilon+\area X,$$
%and so
%$$\vol(\overline{M_0-N})<2\epsilon+\area X.$$
Now since $M=V\cup N\cup\overline{M_0-N}$, we have
$$\begin{aligned}
\vol M&=\vol V +\vol N+\vol(\overline{M_0-N})\cr
&<\epsilon+\epsilon+(\epsilon+\area X)\cr
&=3\epsilon+\area X.
\end{aligned}
$$
Since $\epsilon>0$ was arbitrary, it follows that $\vol M\le\area X$.
\EndProof

%$$M=N\cup C_1\cup\cdots\cup C_n\cup V_1\cup\cdots\cup V_m.$$
%Hence
%\Equation\label{but all the horsemen knew her} 
%\begin{aligned}\vol &M=\vol N+\sum_{j=1}^n\vol %C_j+\sum_{i=1}^m\vol V_i\cr
%&\le\epsilon+\frac12\sum_{j=1}^n\area \partial %C_j+\sum_{i=1}^m\frac\epsilon m\cr
%&\le2\epsilon+\frac12\area (\partial (\overline{M_0-N})).
%\end{aligned}
%\EndEquation
%On the other hand, we have
%$\partial (\overline{M_0-N}))\subset\partial M_0\cup\partial N$, and %hence
%\Equation
%\begin{aligned}
%\area (\partial (\overline{M_0-N}))&\le\area\partial M_0+\area %\partial N\cr
%&<\sum_{i=1}^m\area\partial T_i+(\epsilon+2\area X_0)\cr
%&\le\sum_{i=1}^m\area\frac\epsilon m+\epsilon+2\area X_0\cr
%&=2\area X_0+2\epsilon.
%\EndEquation
%Since $X_0\subset X\cup\partial M_0=

%$$\partial C_1\cup\cdots\cup \partial C_j\subset

%\redcomment{Finish. The following stuff may be of some use:

%This is deduced as follows: Lemmas \ref{oh far out} and \ref{eh golly %you}, with $X=|K|$, imply that if $C$ is any component of $M-f(|K|)$ %and $F_C$ denotes its frontier, we have
%$\vol C\le\frac12\area F_C$. ( We can apply Lemma \ref{eh golly %you} because $\vol M<\infty$. I wrote the following sentence, %whatever it means: ``For the application we need the map to be PL, %which means first `making' $K$ simplicial.'') Summing over the %components of $M-f(|K|)$, we find that
%$$\vol M=\sum_C\vol C\le\frac12\sum_C\area F_C=\area f(|K|),$$
%as required. }
%\EndProof

\begin{proposition}\label{short to bounded}
Let
$(M,\star)$ be a based closed, orientable, hyperbolic $3$-manifold such that $\pi_1(M,\star)$ is generated by two non-commuting elements $x$ and $y$.   Let $\lambda>0$ be given, and suppose that $x$ and $y$ are represented by closed loops of length $<\lambda$ based at $\star$.   Let $W$ be a non-trivial reduced word in two letters such that $W(x,y)=1$.   Then
$$\vol M<(\length (W)-2)\min(\pi, \lambda).$$ (In particular, $M$ has finite volume.)
\end{proposition}

\Proof
If $M$ has infinite volume then Proposition \ref{chained to the altar} implies that the non-commuting elements $x$ and $y$ cannot both be represented by loops of length $<\log3$ based at $\star$. Hence the hypothesis implies that $\vol M<\infty$.

Since $x$ and $y$ do not commute, they are in particular both non-trivial elements of $\pi_1(M,\star)$. Since $\pi_1(M)$ is generated by $x$ and $y$, it is non-abelian and in particular non-cyclic. As the fundamental group of a hyperbolic manifold, it is also torsion-free, and hence $x$ and $y$ have infinite order. 

I will regard $W$ as a word in two abstract generators $\xi$ and $\eta$. Let us write $W=\psi_1\cdots\psi_n$, where each $\psi_j$ has the form $\xi^{\pm1}$ or $\eta^{\pm1}$. We may also identify $W$ with a non-trivial element of the free group $F$ on the generators $x$ and $y$. Since $\pi_1(M)$ is generated by $x$ and $y$ and since $W(x,y)=1$, there is an epimorphism from $F/\langle\langle W\rangle\rangle$ to $\pi_1(M,\star)$ that maps $\xi$ to $x$ and $\eta$ to $y$. Set $n=\length W$. 

If there is a non-trivial reduced word $W'$ in two letters such that $n':=\length W'<n$ and such that $W'(x,y)=1$, and if
$\vol M<(n'-2)\min(\pi, \lambda)$, then in particular
$\vol M<(n-2)\min(\pi, \lambda)$. Hence, arguing inductively, we may assume that $n$ is the minimal length of any non-trivial reduced word which gives a relation between $x$ and $y$.

If $n$ were at most $3$, either $F/\langle\langle W\rangle\rangle$ would be cyclic, which is impossible since $\pi_1(M)$ is non-cyclic, or the image of $\xi$ or $\eta$ in $F/\langle\langle W\rangle\rangle$ would have finite order, which is impossible since $x$ and $y$ have infinite order. Hence $n\ge4$.

The hypothesis implies that there are loops $\alpha,\beta:[0,1]\to M$ based at $\star$ such that $\alpha|(0,1)$ and $\beta|(0,1)$ are open geodesics of length $<\lambda$, and such that $[\alpha]=x$ and $[\beta]=y$. These geodesics are non-constant since $x$ and $y$ are non-trivial elements of $\pi_1(M,\star)$. Let $\calg$ be a graph with one vertex $v$ and two closed edges $a$ and $b$. Fix loops $\kappa _a$ and $\kappa _b$ in $\calg$, based at $v$, such that $\kappa_a|(0,1)$ and $\kappa_b|(0,1)$ are homeomorphisms of $(0,1)$ onto $\inter a$ and $\inter b$ respectively. Define a continuous map $\phi:\calg\to M$ by setting $\phi(v)=\star$, $\phi|\inter a=\alpha\circ \kappa _a^{-1}$ and $\phi|\inter b=\beta\circ \kappa _b^{-1}$.

Let us define a map  $\omega:S^1\to \calg$ as follows. Fix $n$ points on $S^1$, labeled in counterclockwise order as $\zeta_0,\ldots,\zeta_{n-1}$. Set $\zeta_n=\zeta_0$. For $1\le j\le n$ let $A_j$ denote the arc which, when oriented counterclockwise, has $\zeta_{j-1}$ and $\zeta_j$ as its initial and terminal points respectively. Fix a homeomorphism $h_j:A_j\to[0,1]$ which maps $\zeta_{j-1}$ and $\zeta_j$ to $0$ and $1$ respectively. Set $\omega(\zeta_j)=\star$ for $0\le j\le n-1$. For $j=1,\ldots,n$, define $\omega|A_i$ to be $\kappa_a\circ h_j$, $\overline{\kappa_a}\circ h_j$, $\kappa_b\circ h_j$, or $\overline{\kappa_b}\circ h_j$, according to whether
$\psi_j$ is equal to $\xi$, $\xi^{-1}$, $\eta$ or  $\eta^{-1}$, respectively.

Since $W(x,y)=1$, the map $\phi\circ \omega$ is homotopic to a constant. Hence if we define a CW complex $K$ by attaching a $2$-cell to $\calg$ via the attaching map $\omega$, then $\phi$ extends to a map $f:K\to M$. Let $c:D^2\to K$ denote the characteristic map for the $2$-cell of $K$. For each $j$ with $2\le j\le n-2$, let $\sigma_j$ denote the line segment in the Euclidean disk $D^2$ with endpoints $\zeta_0$ and $\zeta_j$. Set $\sigma_1=A_1$ and $\sigma_{n-1}=A_n$. Then for each $j$ with $1\le j\le n$, the topological arc $\sigma_j$ has endpoints $\zeta_0$ and $\zeta_j$. For $j=1,\ldots,n-1$, by precomposing $f\circ c|\sigma_j$ with a homeomorphism from $[0,1]$ to $\sigma_i$ that maps $0$ to $\zeta_0$ and $1$ to $\zeta_j$, we obtain a loop $\nu_j$ based at $\star$. We have $[\nu_j]=W'(x,y)$, where $W':=\psi_1\cdots\psi_j$ is a reduced word of length $j<n$. By our mimimality assumption on $n$, it follows that $\nu_j$ is a homotopically non-trivial loop, and is therefore fixed-endpoint homotopic to a loop whose restriction to $(0,1)$ is a non-constant geodesic. Furthermore, each of the paths $\nu_1|(0,1)$ and $\nu_n|(0,1)$ is a (possibly orientation-reversing) reparametrization of either $\alpha|(0,1)$ or $\beta|(0,1)$, and is therefore a non-constant geodesic. Hence  after modifying $f$ within its homotopy class rel $\calg$ we may assume that $\nu_j|(0,1)$ is a non-constant geodesic for $j=1,\ldots,n-1$.

The set $D^2-(\sigma_1\cup\cdots\cup\sigma_{n-1})$  has $n-2$ components. The closures of these components are topological disks, which we may label as $\tau_1,\ldots,\tau_{n-2}$, where $\partial \tau_j=\sigma_j\cup A_{j+1}\cup\sigma_{j+1}$. 

Let $j$ be any index with $1\le j\le n-2$. The maps $f|\sigma_j$ and $f|\sigma_{j+1}$ are  reparametrizations of the non-constant geodesic paths $\nu_j$ and $\nu_{j+1}$, while $f|A_{j+1}$ may be obtained from one of the non-constant geodesic paths $\alpha$ or $\beta$ by precomposing it with some homeomorphism from $A_{j+1}$ to $[0,1]$. Hence if $\tG_j$ is a lift of the map  $f\circ c| \tau_j$ to $\HH^3$, then $\tG_j(\partial\tau_j)$ is the boundary of a triangle $T_j\subset\HH^3$. Since one of the sides of $T_j$ is a lift of either $\alpha$ or $\beta$, the length of the shortest side of $T_j$ is less than $\lambda$. Hence it follows from Proposition \ref{uneeda} that $\area T_j<\min(\pi,\lambda)$. 

Since $M$ is aspherical, we may arrange after modifying $f$ by a homotopy rel $\calg\cup c(\sigma_2)\cup\cdots\cup c(\sigma_{n-2})$ that for each $j$ with $1\le j\le n-1$, the map $f|\tau_j$ admits a lift to $\HH^3$ which maps $\tau_j$ onto $T_j$. 
Now by \ref{gruosso}  we may fix a smooth triangulation of $M$  with respect to which each of the sets $f(\alpha([0,1])$, $(\beta([0,1])$ and $f\circ c(\tau_j)=p\circ \tG_j(\tau_j)$ ($j=1,\ldots,n-2$), is polyhedral. It follows that the union $f(K)$ of these sets is also polyhedral. We have
$$\area (f\circ c(\tau_j))\le \area T_j<\min(\pi,\lambda).$$
Hence 
\Equation\label{weird love}
\area f(K)\le\sum_{j=1}^{n-2}\area (f\circ c(\tau_j))<(n-2)\min(\pi,\lambda).
\EndEquation

According to the construction of $K$ we have $\pi_1(K)\cong F/\langle\langle W\rangle\rangle$, where, as above, $F$ denotes the free group on the generators $\xi$ and $\eta$. In particular, $\pi_1(K)$ has rank at most $2$. Furthermore, since $W$ is a non-trivial reduced word, $\pi_1(K)$ is not a free group of rank $2$. On the other hand, since $x$ and $y$ generate $\pi_1(M,\star)$, the map $f_\sharp:\pi_1(K,v)\to\pi_1(M,\star)$ is surjective. Since we have observed that $\pi_1(M)$ is non-cyclic, $\pi_1(K)$ is also non-cyclic; hence $\pi_1(K)$ has rank exactly $2$. Thus all the hypotheses of Lemma \ref{oh far out} are seen to hold, and it follows that for each component $C$ of $M-f(K)$, the  image of the inclusion homomorphism $\pi_1(C)\to\pi_1(M)$ is abelian.

We may now apply Lemma \ref{jeweler you failed}, taking $X=f(K)$, to deduce that
\Equation\label{broke in a big way}
\vol M\le\area f(K).
\EndEquation
The required conclusion $\vol M<(n-2)\min(\pi,\lambda)$ follows immediately from (\ref{weird love}) and (\ref{broke in a big way}).
\EndProof

\section{A concrete bound for $\vol M$ when $\mu(M)<(\log3)/2$, and its consequences}\label{concrete section}

\begin{theorem}\label{what, no soap?}
Let $\lambda$ be a positive real number strictly less than $(\log3)/2$, and let $N(\lambda)$ be defined by \ref{weird enough for ya?}. Then every orientable hyperbolic $3$-manifold $M$ with $\mu(M)<\lambda$ we have
$$\vol M< \lambda\cdot (8N(\lambda)-2).$$
\end{theorem}

\Proof
Suppose that $M$ is an orientable hyperbolic $3$-manifold such that $\mu(M)<\lambda$, i.e. such that $\lambda$ is not a Margulis number for $M$. Let us write $M=\HH^3/\Gamma$, where $\Gamma\le\isomplus(\HH^3)$ is discrete and torsion-free. Then there are
non-commuting elements $x$ and $y$ of $\Gamma$ such
that $\max(d(P,x\cdot P),d(P,y\cdot P))<\lambda$ for $i=1,2$.
According to Proposition \ref{first one}, there is a
reduced word $W$ in two letters, with $0<\length  W\le 8N(\lambda)$, such that $W(x,y)=1$.  

Set $\tGamma=\langle x,y\rangle\le\Gamma$, and $\tM=\HH^3/\tGamma$.
Let $\star\in\tM$ denote the image of $P\in\HH^3$ under the quotient map $\HH^3\to\HH^3/\tGamma$. Then $(\HH^3,P)$ is a based covering space of $(\tM,\star)$. Under the natural identification of $\pi_1(\tM,\star)$ with the deck group $\tGamma$, the elements $x$ and $y$ are identified with generators of
$\pi_1(\tM,\star)$ which are represented by closed loops of length $<\lambda$ based at $\star$.   Applying Proposition  \ref{short to bounded}, with $\tM$ playing the role of $M$ in that proposition, we deduce that
$$\vol\tM<\lambda(\length (W)-2)\le \lambda(8N(\lambda)-2).$$ (In particular $\tM$ has finite volume.)
Since $\tM$ covers $M$, it follows that
$$\vol M< \lambda(8N(\lambda)-2).$$
\EndProof

In the following corollary,
$V_0=0.94\ldots$ will denote the volume of the Weeks manifold \cite{weeks}.

\Corollary\label{I can't spell conniptualization}
Let $\lambda$ be a positive real number strictly less than $(\log3)/2$, and let $N(\lambda)$ be defined by \ref{weird enough for ya?}. Then 
for every orientable hyperbolic $3$-manifold $M$, with  $\mu(M)<\lambda$, the group $\pi_1(M)$ has a rank-$2$ subgroup of index at most 
$\lambda\cdot (8N(\lambda)-2)/V_0$.
\end{corollary}

\Proof
Suppose that $M$ is a closed, orientable hyperbolic $3$-manifold 
such that $\mu(M)<\lambda$, i.e. such that $\lambda$ is not a Margulis number for $M$. Write $M_\infty=\HH^3/\Gamma$, where $\Gamma\le\isomplus(\HH^3)$ is discrete and torsion-free. Then by definition there exist a point $P\in\HH^3$ and non-commuting elements $x,y\in\Gamma$ such that
\Equation\label{omg}
\max(d(P,x\cdot P),d(P,y\cdot P))<\lambda.
\EndEquation
Now $\tGamma:=\langle x,y\rangle$ is a non-abelian rank-$2$ subgroup of $\Gamma$.

Since $x$ and $y$ are non-commuting elements of $\tGamma$, it follows from (\ref{omg}) that $\lambda$ is not a Margulis number for $\tM$, i.e. that $\mu(\tM)<\lambda$. Hence by Theorem \ref{what, no soap?} we have $\vol \tM< \lambda\cdot (8N(\lambda)-2)$. (In particular $\tM$ has finite volume) On the other hand, it is shown in \cite{milley} that the Weeks manifold has minimal volume among all orientable hyperbolic $3$-manifolds.  Hence $\vol M\ge V_0$. We therefore have
$$[\Gamma:\tGamma]=\frac{\vol\tM}{\vol M}<
\frac{\lambda\cdot (8N(\lambda)-2)}
{V_0}.$$
It follows that $\pi_1(M)\cong\Gamma$ has a rank-$2$ subgroup of index at most 
$\lambda\cdot (8N(\lambda)-2)/
{V_0}$.
\EndProof

\Corollary\label{talk like a putz day}
Let $\lambda$ be a positive real number strictly less than $(\log3)/2$. Then for every hyperbolic $3$-manifold $M$ with $\mu(M)<\lambda$, we have
$$\rank\pi_1(M)\le2+\log_2(\lambda\cdot (8N(\lambda)-2)/V_0).$$
\end{corollary}

\Proof
Let $M$ be an orientable hyperbolic $3$-manifold with
$\mu(M)<\lambda$. According to Corollary \ref{I can't spell conniptualization}, $\pi_1(M)$ has a rank-$2$ subgroup $X$ of index at most 
$\lambda\cdot (8N(\lambda)-2)/V_0$. According to Proposition \ref{you peeked}, we have
$$\rank\pi_1(M)\le \rank X+\log_2[\pi_1(M):X]\le2+
\log_2\bigg(\frac{\lambda\cdot (8N(\lambda)-2)}{V_0}\bigg).$$
\EndProof

I will conclude with three corollaries which follow immediately upon combining the earlier results of this section with the estimate for $N(\lambda)$ given by Proposition \ref{nestimate} and Meyerhoff's lower bound \cite{meyerhoff} of $0.104$ for $\mu_+(3)$. The first, Corollary \ref{boingo cuckoo}, was stated in the introduction as Theorem B, and the other two were presented as corollaries to Theorem B (although I will derive them formally from Corollaries \ref{I can't spell conniptualization} and \ref{talk like a  putz day} above.)

\Corollary\label{boingo cuckoo} 
Let $\lambda$ be a positive real number strictly less than $(\log3)/2$. Then for every orientable hyperbolic $3$-manifold $M$ with $\mu(M)<\lambda$ we have
$$\vol M< \lambda\bigg(6+\frac{880}{\log3-2\lambda}\log{1\over\log3-2\lambda}\bigg).$$
\EndCorollary

\Proof
Since $\lambda>\mu(M)$, we have in particular that $\lambda>\mu_+(3)>0.1$.
The assertion now follows from Theorem \ref{what, no soap?} and Proposition \ref{nestimate}.
\EndProof

\Corollary\label{more cuckoo}
Let $\lambda$ be a positive real number strictly less than $(\log3)/2$. Then 
for every orientable hyperbolic $3$-manifold $M$ with $\mu(M)<\lambda$, the group $\pi_1(M)$ has a rank-$2$ subgroup of index at most 
$$\frac{\lambda}{V_0}\bigg(6+\frac{880}{\log3-2\lambda}\log{1\over\log3-2\lambda}\bigg).$$
\EndCorollary

\Proof
Since $\lambda>\mu(M)$, we have in particular that $\lambda>\mu_+(3)>0.1$.
The assertion now follows from Corollary \ref{I can't spell conniptualization} and Proposition \ref{nestimate}.
\EndProof

\Corollary\label{mostly moxie}
Let $\lambda$ be a positive real number strictly less than $(\log3)/2$. Then for every hyperbolic $3$-manifold $M$ with 
$\mu(M)<\lambda$, we have
$$\rank\pi_1(M)\le2+\log_2\bigg(\frac{\lambda}{V_0}\bigg(6+\frac{880}{\log3-2\lambda}\log{1\over\log3-2\lambda}\bigg)\bigg).$$
\EndCorollary

\Proof
Since $\lambda>\mu(M)$, we have in particular that $\lambda>\mu_+(3)>0.1$.
The assertion now follows from Corollary \ref{talk like a putz day} and Proposition \ref{nestimate}.
\EndProof

\bibliographystyle{plain}
\bibliography{bounds}

\end{document}